\documentclass[12pt]{amsart}
\usepackage{amssymb,amsmath,bbm, amsthm}
\usepackage{graphicx} 
\usepackage{tikz}

\usepackage{booktabs,array,tabularx}
\usepackage{todonotes}
\usetikzlibrary{arrows.meta}

\newcommand{\Xb}{\mathbf{X}}
\newcommand{\xb}{\mathbf{x}}
\newcommand{\yb}{\mathbf{y}}
\newcommand{\zb}{\mathbf{z}}
\newcommand{\cb}{\mathbf{c}}
\newcommand{\ab}{\mathbf{a}}
\newcommand{\kb}{\mathbf{k}}
\newcommand{\mb}{\mathbf{m}}
\newcommand{\jb}{\mathbf{j}}
\newcommand{\jk}{\mathfrak{j}}

\newcommand{\pr}{{\rm pr}}
\newcommand{\im}{{\rm im}}
\newcommand{\pd}{{\rm pd}}
\newcommand{\sign}{{\rm sign}}
\newcommand{\cl}{\mathbf{cl}}
\newcommand{\rank}{{\rm rank}}
\newcommand{\supp}{{\rm supp}}

\newcommand{\NN}{\mathbb{N}}
\newcommand{\ZZ}{\mathbb{Z}}
\newcommand{\FF}{\mathbb{F}}
\newcommand{\1}{\mathbbm{1}}
\newcommand{\TT}{\mathbb{T}}
\newcommand{\Sc}{\mathcal{S}}
\newcommand{\Cc}{\mathcal{C}}
\newcommand{\Fc}{\mathcal{F}}
\newcommand{\Hc}{\mathcal{H}}
\newcommand{\Nc}{\mathcal{N}}
\newcommand{\Pc}{\mathcal{P}}
\newcommand{\HN}{H\!N}

\newcommand{\lcm}{{\rm lcm}}

\newtheorem{theorem}{Theorem}[section]
\newtheorem{definition}{Definition}[section]
\newtheorem{corollary}{Corollary}[section]
\newtheorem{proposition}{Proposition}[section]
\newtheorem{example}{Example}[section]
\newtheorem{remark}{Remark}[section]

\makeatletter

\title[Algebraic and topological tools for domination theory]{Algebraic and topological tools for the analysis of domination theory for multi-state systems}

 \author{Patricia Pascual-Ortigosa}
 \address{Universidad de La Rioja (Spain)}
 \email{papasco@unirioja.es}
 \author{Eduardo Sáenz-de-Cabezón}
 \email{esaenz-d@unirioja.es}

 \date{\today}

\begin{document}

\maketitle

 \begin{abstract}

 The signed domination function of a coherent system is the M\"obius
 transform of its structure function. It gives the coefficients of the
 reliability polynomial and measures the intrinsic difficulty of exact
 reliability computations. We study it for binary and multi-state systems by
 means of the monomial ideals associated with its performance levels.
 
 Our main result states that, at every level and every multidegree, the
 signed domination function equals the multigraded Euler characteristic of
 the corresponding reliability ideal, that is, the alternating sum of its
 multigraded Betti numbers, which is the corresponding coefficient of its
 multigraded Hilbert numerator. Equivalently, it equals minus the reduced
 Euler characteristic of the Koszul simplicial complex of the ideal at that
 multidegree.
 
 This yields two computational methods: one based on Koszul complexes,
 requiring at most $2^n$ evaluations of the structure function for a system
 with $n$ components, and one based on Mayer--Vietoris trees, whose size is at
 most $2^r$ for a system with $r$ minimal paths. It also shows that computing
 the domination of a system is equivalent to computing the Euler
 characteristic of a simplicial complex given by its facets, a \#P-complete
 problem. Polarization and depolarization relate the dominations of
 multi-state and binary systems.
 
 We distinguish the domination at the top state from the level domination,
 and give a sufficient condition for zero domination. As applications, we
 compute the signed domination of component-based dominant systems,
 characterize algebraically dominant multi-state systems, and determine the dominations of binary and multi-state $k$-out-of-$n$
 and consecutive $k$-out-of-$n$ systems, in simple and generalized versions.
 
 \end{abstract}



\section{Introduction}
The structure function of binary coherent systems can be studied by means of {\em domination theory} \cite{SP78, B82, SC83}. It can be expressed as products of states of the components of the system, with coefficients given by the so called {\em signed domination function}. Formally, the signed domination function is the M\"obius transform of the structure function of the system, and is usually described and computed in terms of odd and even formations of elements of the boolean lattice (see Definition \ref{def:signedDomination}). It was introduced by Satyanarayana and Prabhakar \cite{SP78} in their topological formula for computing network reliability, see also \cite{AB84}. The most influential application of domination theory is as a measure of the intrinsic difficulty of exact reliability computation. In \cite{SC83} it is shown that the absolute value of the signed domination of an undirected $K$-terminal network is the number of irreducible leaves in an optimal factoring strategy. See \cite{SW85,W86,H89} for further algorithmic developments and a general framework for coherent systems. Domination theory has also been linked to signatures of coherent systems, see \cite{B01,MM11}, and to importance measures, in relation to simple games \cite{R90}. The theory of signed domination functions has been recently extended in \cite{H25} to multi-state systems in a systematic way.

On the other hand, the theory of monomial ideals provides useful tools to study coherent systems (both binary and multi-state), see e.g. \cite{SW09,SW15}. A main insight of this algebraic approach is the use of the Hilbert series of the so called reliability ideal of a system, to obtain compact inclusion-exclusion formulas for the reliability of the system, both in the binary and the multi-state case \cite{MPSW20}.

In the present paper we connect these two theories, and give an algebraic description of domination theory by means of algebraic and topological tools based on monomial ideals and simplicial complexes related to them. The main result of the paper is the connection between domination theory and monomial ideals by the following identity (see Theorem \ref{th:dominationBetti}) for each level $j$ of performance of the system,
\begin{equation}\label{eq:identity}
    \delta_{\Sc,j}(\mu)=\sum_i (-1)^i\beta_{i,\mu}(I_j(\Sc))=\chi_\mu(I_j(\Sc)),
\end{equation}

where $\Sc$ is a coherent multi-state system, $\delta_{\Sc,j}$ is its signed domination function at level $j$, $\mu\in\NN^n$ is a state vector of $\Sc$,  $I_j(\Sc)$ is the system's $j$-reliability ideal, $\beta_{i,\mu}(I_j(\Sc))$ is the $i$-th Betti number of $I_j(\Sc)$ at multi-degree $\mu$, and $\chi_\mu(I_j(\Sc))$ is the multi-graded Euler characteristic of $I_j(\Sc)$ at multi-degree $\mu$.

The Betti numbers of a monomial ideal are its main homological invariants, which are used to describe the structure of the ideal. The combinatorial nature of these ideals allow the use of tools from other areas, in particular simplicial complexes and their Euler characteristic, which are the main object in combinatorial topology \cite{S96,MS05}. We explore some of these tools in order to take advantage of the identity (\ref{eq:identity}) for the study and computation of signed domination functions of coherent systems.

The outline of the paper is the following. In Section \ref{sec:preliminaries} we give the necessary background on coherent systems, monomial ideals and the relations between them. We also describe the basics of domination theory. Section \ref{sec:combinatorics} describes several combinatorial objects related to monomial ideals, including the $\lcm$-lattice, Taylor resolution and Stanley-Reisner and Koszul simplicial complexes. Section \ref{sec:dominationEuler} relates the signed domination function of system $\Sc$ to the Euler characteristic of certain chain and simplicial complexes related to the reliability ideal of the system $\Sc$. Section \ref{sec:computationMVT} computes the Euler characteristic, and Section \ref{sec:polarization} studies the connection between binary and multi-state systems by means of polarization and depolarization. In Section \ref{sec:examples} we apply the introduced techniques to compute the domination of several binary and multi-state coherent systems. We define and describe dominant systems, and we cover binary and multi-state $k$-out-of-$n$ and consecutive linear $k$-out-of-$n$ systems.

\section{Preliminaries}\label{sec:preliminaries}
\subsection{Coherent systems}\label{sec:coherent}
A multi-state system $\Sc=(C,\phi),$ consists of a set $C=\{\bar{c}_1,\dots, \bar{c}_n\}$ of components and a structure function $\phi$. Each component $\bar{c}_i$ can be in a finite number of states $S_i=\{0,1,\dots,m_i\}$, that represent increasing levels of performance of the component, where $0$ represents total failure, and $m_i$ represents full performance. The space of states of the components is denoted as $\Cc=S_1\times\cdots\times S_n$. A system state vector is an element $\cb=(c_1,\dots,c_n)$ of $S_1\times\cdots\times S_n$, where each $c_i\in S_i$ represents the state of component $\bar{c}_i$. Similarly, the system itself can be in a finite number of increasing states of performance $S=\{0, 1,\dots,M\}$. If $m_i=1$ for all $i$ and $M=1$ then we say that the system is binary. If $M>1$ or there is some $i$ for which $m_i>1$ we say that the system is multi-state. The structure function $\phi$ expresses the state of the system as a function of the states of the components, i.e. $\phi = \phi(\cb): \Cc \rightarrow S$. The system $\Sc$ is called {\em monotone} if $\phi$ is non-decreasing on each component. $\Sc$ is called {\em coherent for level $j\in\{0,\dots,M\}$} if it is monotone and each component is relevant for level $j$. If all components are relevant for all levels, the system is said to be a {\em coherent} system.

A state vector $\cb\in\Cc$ is a $j$-level {\em path vector} for system, $\Sc=(C,\phi)$, if $\phi(\cb)\geq j$. We say that $\cb$ is a minimal $j$-level path vector of $\Sc$ if $\phi(\cb)\geq j$ and $\phi(\cb')<j$ for all state vectors $\cb'\in\Cc$ such that $\cb'<\cb$. We denote by $\Fc_j(\Sc)$ and $\overline{\Fc}_j{(\Sc})$ the set of $j$-path vectors and minimal $j$-path vectors of $\Sc$ respectively.

To each component $\bar{c}_i$ and level $j\in\{0,\dots,m_i\}$ we assign the probability $\pr_{i,j}:=\pr(c_i \geq j)$ that component $\bar{c}_i$ is in a state better than or equal to $j$. The probability that the state vector of the variables of the system is bigger than a certain vector $\jb=(j_1,\dots,j_n)$ is $\pr(\jb)=\pr((c_1\geq j_1)\land\cdots\land(c_n\geq j_n))$. If all probabilities are independent, then $\pr(\jb)=\prod_{i=1}^n\pr_{i,j_i}$. The {\em $j$-reliability of system $\Sc$}, denoted as $Rel_j(\Sc)$, or just $Rel(\Sc)$ in the binary case, is the probability that the system is at least in level $j\in\{0,\dots,M\}$. The reliability of the system is computed using the structure function $\phi$ of $\Sc$. For a complete introduction of the theory and methods of reliability theory, see for instance \cite{KZ03,TB17}.

\subsection{Domination theory}
The main way to compute system reliability is to use the structure of the set of paths of the system. The problem of reliability computation is \#P-hard \cite{B86}, hence algorithms and tools that exploit the structure of the problem through mathematical tools are of paramount importance. One of these tools is {\em domination theory}, that plays an important role in e.g. the computation of structure functions, the factoring algorithm and undirected network systems \cite{SP78,SC83,SW85, W86}. The main concepts in this theory are the following (the notations in this section follow those in \cite{H25}).

Let $x,y\in\NN$ and denote $x\lor y:= \max\{x,y\}$. If $\xb=(x_1,\dots,x_n)$ and $\yb=(y_1,\dots,y_n)$ are vectors in $\NN^n$, then we denote $\xb\lor\yb:=(x_1\lor y_1,\dots,x_n\lor y_n)$. We say that $\xb\leq\yb$ if $x_i\leq y_i$ for all $i=1,\dots,n$. If $\xb\nleq\yb$ and $\yb\nleq\xb$ then we say that $\xb$ and $\yb$ are incomparable. 

Let $\Xb=\left\lbrace\xb_1,\dots,\xb_r\right\rbrace$ be a set of incomparable vectors in $\NN^n$. We say that the {\em maximum-closure} of $\Xb$, denoted $\cl(\Xb)$, is the minimal set of vectors in $\NN^n$ that contains $\xb_1,\dots \xb_r$ and that is closed under the $\lor$ operation, i.e. $\xb,\yb\in\cl(\Xb)\implies\xb\lor\yb\in\cl(\Xb)$. The set $\cl(\Xb)$ is a poset with respect to $\leq$, where $\xb_i, i=1,\dots,r$ are the minimal elements and $\xb_1\lor\xb_2\lor\dots\lor\xb_r$ is the (unique) maximal element.

\begin{definition}\label{def:signedDomination}
Let $\xb\in \cl(\Xb)$, then each set $\{\xb_{i_1},\dots,\xb_{i_f}\}$ of minimal elements of $\cl(\Xb)$ such that $\xb=\xb_{i_1}\lor\dots\lor\xb_{i_f}$ is called a {\em formation} of $\xb$. If $f$ is odd, we say that this formation is {\em odd}, otherwise we say it is an {\em even formation}.

The {\em signed domination function} of $\cl(\Xb)$, $\delta(\xb)$, is defined as
\begin{equation}\label{eq:signedDomination}
\delta(\xb)=\vert \mbox{odd formations of }\xb \vert-\vert \mbox{even formations of }\xb \vert, \mbox{ for }\xb\in\cl(\Xb),
\end{equation}
and we set $\delta(\xb)=0$ for $\xb\notin\cl(\Xb)$.
\end{definition}

Observe that the signed domination function is invariant under poset isomorphisms.

When we choose as $\Xb$ the set of minimal $j$-paths of a multi-state system $\Sc$, $\overline{\Fc}_j(\Sc)$, we use the following notation \[\delta_{\Sc,j}(\xb)=\vert \mbox{odd formations of }\xb \vert-\vert \mbox{even formations of }\xb \vert, \mbox{ for }\xb\in\cl(\overline{\Fc}_j(\Sc)).\]
When the system $\Sc$ is clear, we simply denote $\delta_{\Sc,j}(\xb)$ by $\delta_{j}(\xb)$ to ease the notation.

Let $\Sc=(C,\phi)$ be a coherent system. The signed domination function has a strong relation with the structure function of the system. 
\begin{definition}
    Let $\phi_j$ be the binary function $\phi_j(\cb)=I(\phi(\cb)\geq j)$, i.e. the indicator function of the $j$-paths of $\Sc$. The functions $\phi_j$ for $j\in\{1,\dots,M\}$ are called the \emph{$j$-level structure functions of $\Sc$}.
\end{definition}

\begin{theorem}[\cite{H25}, Th. 3.2]
Let $\Sc=(C,\phi)$ be a multi-state coherent system and let $\overline{\Fc}_j(\Sc)$ be the set of its minimal $j$-path vectors. The $j$-level structure function $\phi_j$ can be written as 
\[
\phi_j(\yb)=\sum_{\xb\in\cl(\overline{\Fc}_j(\Sc))}\delta_j(\xb)I(\yb\geq\xb),\;\forall\, \yb\in\Cc,
\]
where $\delta_j$ is the signed domination function of $\cl(\overline{\Fc}_j(\Sc))$, and $I(\yb\geq \xb)=1$ if $\yb\geq \xb$, and $0$ otherwise.
\end{theorem}

\begin {definition}
Let $\Sc=(C,\phi)$ be a multi-state system and $\Cc$ its state space.
The signed domination of the maximal element $\mb$ of $\Cc$ is called the {\em signed domination} of the structure function $\phi_j$, and denoted as 
\[
d(\phi_j)=\delta_j(\mb)=d_j(\Sc).
\]

The signed domination of the maximal element of $\cl(\overline{\Fc}_j(\Sc))$ is called the {\em level domination} of the structure function $\phi_j$, and is denoted as 

\[
\overline{d}(\phi_j)=\delta_j\left(\bigvee{\overline{\Fc}_j(\Sc)}\right)=\overline{d}_j(\Sc).
\]
\end{definition}

\subsection{Monomial ideals and Algebraic Reliability}\label{sec:monomial}

Let $\Sc=(C,\phi)$ be a coherent system with $n$ components. Consider the polynomial ring $R=\kb[x_1,\dots,x_n]$ in $n$ variables, where $\kb$ is a field (in this paper we consider $\kb=\mathbb{R}$). To each vector state $\cb=(c_1,\dots,c_n)$ in $\Cc$ we assign the monomial $\xb^\cb=x_1^{c_1}\cdots x_n^{c_n}$. The monomial ideal generated by the monomials corresponding to the minimal $j$-paths of $\Sc$ is called the {\em $j$-reliability ideal} of $\Sc$ and is denoted by $I_j(\Sc)$, i.e.
\[
I_j(\Sc)=\left\langle \xb^\cb \vert \cb\in \overline{\Fc}_j(\Sc)\right\rangle.
\]

This equivalence between coherent systems and monomial ideals was used in a series of papers, see e.g. \cite{GW04,SW09,SW10,SW15,BPS21} to develop an algebraic approach to binary and multi-state system reliability. In order to compute the $j$-reliability of a system $\Sc$ for a certain level $j$, the algebraic approach needs to efficiently list the monomials corresponding to all $j$-paths of $\Sc$. This counting is given by the Hilbert function and Hilbert series of the corresponding ideal. We give here a brief account of these concepts; for a detailed treatment, see \cite{HH11,V15}.

Let $I\subseteq\kb[x_1,\dots,x_n]$ be a monomial ideal and $G(I)=\{\xb^{\mu_1},\dots,\xb^{\mu_r}\}$ its minimal generating set. The Hilbert function $H_I$ of $I$ is an integer valued function that assigns to each integer $d$ the dimension of the $\kb$-vector space generated by the monomials in $I$ of degree $d$, i.e. the number of monomials in $I$ of degree $d$. The Hilbert series of $I$, denoted by $\Hc_I$ is the generating function of $H_I$. By the celebrated Hilbert-Serre theorem, see \cite{AM69}, Theorem 11.1, the Hilbert series is a rational function. 
Both the Hilbert function and series are multi-graded, hence we can have a multi-graded version of them. In the case of the Hilbert series, we can express it in the form
\[
\Hc_I=\frac{\HN_I(x_1,\dots,x_n)}{\prod_{i=1}^n (1-x_i)},
\]
where $\HN_I(x_1,\dots,x_n)$ is a polynomial in $x_1,\dots,x_n$ that we will call {\em (multi-graded) Hilbert numerator of $I$}. The Hilbert numerator of $I_j(\Sc)$ gives us the $j$-reliability polynomial of the system $\Sc$, which is the core idea in Algebraic Reliability, see the details in e.g \cite{SW09,SW10,SW15}. Therefore, an efficient computation of Hilbert numerators of monomial ideals is a valuable tool in the algebraic approach to system reliability. There are plenty of different ways to compute $\HN_I$. Some of them, involving rich combinatorial structures, will be presented in Section \ref{sec:combinatorics}, which will be later used in Section \ref{sec:dominationEuler} to establish a connection between the algebra of monomial ideals and domination theory.

Another important object that describes the structure of monomial ideals is a free resolution \cite{CLO05,HH11,V15}. Let $I\subseteq R=\kb[x_1,\dots,x_n]$ be a monomial ideal, a {\em free resolution} of $I$ is a sequence of free modules and morphisms between them (called {\em differentials}) of the form
\begin{equation}\label{eq:res}
0\longrightarrow R^{r_l}\stackrel{\partial_l}{\longrightarrow}R^{r_{l-1}}\stackrel{\partial_{l-1}}{\longrightarrow}\cdots \longrightarrow R^{r_1}\stackrel{\partial_1}{\longrightarrow}R^{r_0}\longrightarrow I\longrightarrow 0,
\end{equation}
where $l$ is the length of the resolution, $r_0$ is the number of generators of $I$ and $\ker(\partial_i)=\im(\partial_{i+1})$ for all $i$. Since $I$ is a multi-graded module, each of the modules $R^{r_i}$ is of the form $\oplus_{j=1}^{j=r_i}R(-\mu_{i,j})$ for $\mu_{i,j}\in\NN^n$, i.e. each of the modules is also multi-graded. Due to the additivity of the Hilbert function and series on exact sequences and the rank-nulity theorem, for any free resolution of the form (\ref{eq:res}) we have that
\begin{equation}\label{eq:resHilbert}
 \HN_I=\sum_{i=0}^{i=l}(-1)^i\sum_{j=1}^{j=r_i}\xb^{\mu_{i,j}},
\end{equation}
i.e. free resolutions provide a way to compute the Hilbert numerator of $I$. This means that free resolutions are tools for computing the reliability polynomial and bounds for binary and multi-state coherent systems \cite{GW04,SW09}. Among all possible free resolutions of a monomial ideal $I$, there is a unique minimal one (up to isomorphism). In the case of a minimal resolution, the ranks $r_i$ receive the name of Betti numbers and are denoted by $\beta_i(I)$ since they are invariants of the ideal. When taking into account the multi-graded structure, we consider the multi-graded Betti numbers of $I$, and denote them $\beta_{i,\mu}(I)$. The length of the minimal resolution is called the {\em projective dimension} of $I$ and denoted as $\pd(I)$. The alternating sum of the Betti numbers is called the Euler characteristic of the ideal.
\begin{definition}\label{def:EulerCharIdeals}
    Let $I\subseteq \kb[x_1,\dots, x_n]$ be a monomial ideal. The {\em Euler characteristic} of $I$ is defined as
    \[
    \chi(I):=\sum_{i=0}^n(-1)^i\beta_i(I).
    \]
Taking into account the multi-graded structure, we define the {\em multi-graded Euler characteristic} of $I$
 at multi-degree $\mu$ as
 \[
 \chi_\mu(I):=\sum_{i=0}^n(-1)^i\beta_{i,\mu}(I).
 \]
 \end{definition}

 When $I=I_j(\Sc)$ is the $j$-reliability ideal of a multi-state coherent system $\Sc$, then the numerator of the Hilbert series of $I_j(\Sc)$ gives the $j$-reliability of $\Sc=(\Cc,\phi)$ by the substitution of each monomial $\xb^\mu$ in $\HN(I_j(\Sc))$ by the probability associated to its exponent; i.e. given $\xb^\mu=x_1^{\mu_1}\cdots x_n^{\mu_n}$,
\[
\pr(\xb^\mu)=\pr(\mu)=\pr((c_1\geq \mu_1)\land\cdots\land(c_n\geq \mu_n)),
\]
which, as we said before, in the case of independent probabilities can be obtained as
\[
\pr(\xb^\mu)=\prod_{i=1}^n \pr_{i,\mu_i}.
\]

Using this substitution, we obtain the $j$-reliability of system $\Sc$ as a probabilistic form of expression (\ref{eq:resHilbert})
\begin{equation}\label{eq:resHilbertReliability}
 Rel_j(\Sc)=\sum_{i=0}^{i=l}(-1)^i\sum_{j=1}^{j=r_i}\pr(\mu_{i,j}),
\end{equation}
where the $\mu_{i,j}$ correspond to the multi-degrees in a resolution of $I_j(\Sc)$.

Truncating (\ref{eq:resHilbertReliability}) at even sumands, we obtain upper bounds for the $j$-reliability of $\Sc$, and truncating at odd summands we obtain lower bounds. These are a generalization of the so called {\em Bonferroni bounds} in the literature, see e.g. \cite{D03,NW05}. Among the bounds provided by free resolutions for the reliability of coherent systems, the ones given by the minimal free resolution are the tightest ones, and the reliability polynomial is the most compact within this class, providing significant improvements in terms of reducing redundancy \cite{SW09,SW10,IMPSW24}.

\begin{example}
Consider a binary system $\Sc$ with $n$ components that works whenever at least $k$ of them are working. These systems are called $k$-out-of-$n$:G systems (G is for {\em good}), see Section \ref{sec:k-out-of-n} below. Let's take $n=5,k=3$. The reliability ideal of $\Sc$ is given by
\begin{align*}
I_\Sc=\langle &x_1x_2x_3,x_1x_2x_4,x_1x_2x_5,x_1x_3x_4,x_1x_3x_5,x_1x_4x_5,\\
&x_2x_3x_4,x_2x_3x_5,x_2x_4x_5,x_3x_4x_5\rangle,
\end{align*}
whose minimal generators correspond to the minimal paths of $\Sc$. The numerator of the Hilbert series of $I_\Sc$, as given by its minimal free resolution, is
\begin{align*}
HN_{I_\Sc}&=x_1x_2x_3+x_1x_2x_4+x_1x_2x_5+x_1x_3x_4+x_1x_3x_5+x_1x_4x_5+x_2x_3x_4\\
&+x_2x_3x_5+x_2x_4x_5+x_3x_4x_5-3(x_1x_2x_3x_4+x_1x_2x_3x_5+x_1x_2x_4x_5\\
&+x_1x_3x_4x_5+x_2x_3x_4x_5)+6(x_1x_2x_3x_4x_5).
\end{align*}
Then, if $pr_{i,1}=p$ for all $i$, we have that
\[
Rel(\Sc)=10p^3-15p^4+6p^5.
\]
\end{example}

\section{Combinatorial objects related to monomial ideals}\label{sec:combinatorics}
Monomial ideals lie in the interplay between algebra and combinatorics. Their {\em algebraic} nature implies that we have at hand all the tools and concepts of commutative algebra to study them, such as dimension theory, homological methods, etc. These include Hilbert function and series and free resolutions among others. On the other hand, the {\em combinatorial} nature of monomial ideals allows the use of multiple other combinatorial objects for their study. In this section we review some of these concepts with an emphasis on those that are useful for the study of domination theory.

\subsection{Simplicial complexes: Stanley-Reisner correspondence and the simplicial Koszul complex.}
An abstract simplicial complex $\Delta$ on the set of vertices $[n]=\{1,\dots,n\}$ is a collection of subsets of $[n]$ such that if $\sigma\in\Delta$ and $\tau\subseteq\sigma$ then $\tau\in \Delta$, i.e. $\Delta$ is closed under taking inclusions. 
 The elements of $\Delta$ are called faces of $\Delta$ and the maximal faces with respect to inclusion are called {\em facets}. Observe that a simplicial complex is completely determined by its set of facets. A face $\sigma$ has dimension $d$, denoted $\dim(\sigma)=d$, if $d =\vert\sigma\vert-1$. The dimension of the empty set is $-1$. The dimension $\dim(\Delta)$ of $\Delta$ is defined as $\max\{\dim(\sigma)\vert \sigma\in\Delta\}$. A set $\sigma\notin \Delta$ is called a non-face of $\Delta$; the set of non-faces of $\Delta$ is denoted $\Nc(\Delta)$. The $f$-vector of a simplicial complex is the vector $(f_{-1},f_0,f_1,\dots,f_d)$ where $f_i$ is the number of faces of dimension $i$ in $\Delta$. The {\em reduced Euler characteristic} of $\Delta$, $\tilde \chi(\Delta)$ is a fundamental invariant of $\Delta$, and is defined as \footnote{In this context, we use the {\em reduced} Euler characteristic, starting the definition at dimension $-1$ instead of $0$. In a simplicial complex, the $-1$-dimensional simplex is defined as the empty set ($\emptyset$). Because every non-empty simplicial complex contains exactly one empty set, $f_{-1} = 1$. Hence, for every non-empty simplicial complex $\Delta$, $\tilde{\chi}(\Delta)=\chi(\Delta)-1$.} 
 \[
 \tilde{\chi}(\Delta)=\sum_{i=-1}^d(-1)^if_i.
 \]
 Simplicial complexes are the main object in combinatorial algebraic topology; for general definitions and the basic theory of simplicial complexes see e.g. \cite{H04}.

Given a simplicial complex of dimension $d$ and a field $\kb$, we can define a chain complex of $\kb$-vector spaces $C_i$, $i=0,\dots,d$, where each $C_i$ is generated by generators $e^i_\sigma$ corresponding to the faces $\sigma$ of dimension $i$ of $\Delta$. The chain complex has the form
\[
0\longrightarrow C_d\stackrel{\partial_d}{\longrightarrow} C_{d-1}\stackrel{\partial_{d-1}}{\longrightarrow} \cdots C_2\stackrel{\partial_2}{\longrightarrow} C_1 \stackrel{\partial_1}{\longrightarrow} C_0\longrightarrow \kb\longrightarrow 0,
\]
where each operator $\partial_i$ is called the $i$'th boundary operator, and is given by
\[
\partial_i(e^i_\sigma)=\sum_{s\in\sigma}(-1)^{\sign(s,\sigma)}e^{i-1}_{\sigma\setminus s},
\]
where ${\sign(s,\sigma)}$ is the position of index $s$ in $\sigma$. It is easy to see that $\im(\partial_{i+1})\subseteq \ker(\partial_i)$ for all $i$; the quotient $\tilde{H}_i(\Delta,\kb):=\ker(\partial_i)/\im(\partial_{i+1})$ is the $i$-th (reduced) homology module of $\Delta$, and its dimension $\dim_\kb(\tilde{H}_i(\Delta,\kb))$ is called the $i$-th Betti number of $\Delta$ with coefficients in the field $\kb$. These are usually denoted by $\beta_i(\Delta,\kb)=\dim_\kb \tilde{H}_i(\Delta,\kb)$. Note that although they are related (see Section \ref{sec:Hochster}), the Betti numbers of monomial ideals and those of simplicial complexes are two different notions, defined in different ways. By the Euler-Poincaré theorem \cite{P93}, we have that
\begin{equation}
    \tilde{\chi}(\Delta)=\sum_{i=-1}^d (-1)^{i}\beta_i(\Delta,\kb).
\end{equation}
There are two main ways to associate monomial ideals and simplicial complexes. Historically, the most important one is the Stanley-Reisner correspondence, which is a one-to-one correspondence between simplicial complexes and squarefree monomial ideals \cite{S96} given as follows.
\begin{definition}
    Let $\Delta$ be a simplicial complex, and for each $\sigma\subset[n]$ let $\xb_\sigma=\prod_{i\in\sigma}x_i$. The ideal
    \[
    I_\Delta=\langle \xb_\sigma\vert\sigma\in\Nc(\Delta)\rangle
    \]
 is called the {\em Stanley-Reisner ideal of $\Delta$}. The ring $R/I_\Delta$ is called the Stanley-Reisner ring of $\Delta$.
\end{definition}
This correspondence was used by Stanley to prove the Upper Bound Conjecture for spheres \cite{S75, S17}, a seminal result that opened a fruitful interaction between combinatorial commutative algebra and algebraic topology.

A second construction relating monomial ideals and simplicial complexes (in fact a generalization of the previous one) is the {\em Koszul simplicial complex}, see \cite{MS05}, which is useful for the fast computation of the homology of simplicial complexes using monomial ideals \cite{BHS19}.

\begin{definition}\label{def:KoszulComplex}
Let $I\subseteq R=\kb[x_1,\dots,x_n]$ be a monomial ideal and $\mu\in\NN^n$ a multi-index. Consider the simplicial complex $\Delta K(I,\mu)$ given by $\sigma\in\Delta K(I,\mu)$ iff $\frac{\xb^\mu}{\xb^\sigma}\in I$ with $\sigma\subseteq\supp(\mu)$. The complex $\Delta K(I,\mu)$ is called the {\em Koszul simplicial complex of $I$ at $\mu$}.
\end{definition}

Conversely, to any simplicial complex $\Delta$ one can associate its Koszul ideal as follows.
\begin{definition}\label{def:KoszulIdeal}
Let $\Delta$ be a simplicial complex in $n$ vertices, the {\em Koszul ideal} of $\Delta$ is the square-free monomial ideal
\[
I^{\Delta} = \langle x_{\overline{\sigma}}\ \vert\, \sigma \mbox{ is a facet of } \Delta\rangle,
\]
where $\overline{\sigma}$ is the complement of $\sigma$ in $[n]$.
\end{definition}

\subsection{The lcm-lattice}\label{sec:lcm}
Let $I$ be a monomial ideal and $G(I)=\{\xb^{\mu_1},\dots,\xb^{\mu_r}\}$ its minimal monomial generating set. For each $\sigma\subset\{\mu_1,\dots,\mu_r\}$ we denote by $\xb^{\mu_\sigma}$ the least common multiple of those generators of $I$ that are in $\sigma$, i.e. $\mu_\sigma=\lcm(\{\xb^{\mu_i}\vert \mu_i\in\sigma\})$. By an abuse of notation, we also consider that given a subset $\sigma\subseteq[r]$, $\xb^{\mu_\sigma}=\lcm(\{\xb^{\mu_i}\vert i\in\sigma\})$. The monomials $\xb^{\mu_\sigma}$ for each $\sigma\subset \{\mu_1,\dots,\mu_r\}$ form a poset with respect to divisibility which is a lattice, called the {\em $\lcm$-lattice} of $I$, denoted by $L_I$. The $\lcm$-lattice of a monomial ideal was introduced in \cite{GPW99} and is one of the main tools for the study of certain types of resolutions of monomial ideals, see e.g. \cite{P02,C10,M13,FMS16}. Of particular importance to our goals is the fact that $\beta_{i,\mu}(I)=0$ for all $\xb^{\mu}\notin L_I$, meaning in particular that the set of multi-degrees at which a monomial ideal might have non-zero Betti numbers is finite.

Another crucial fact is that the Betti numbers of $I$ are determined by the poset homology of $L_I$, see Section \ref{sec:Hochster}. Associated to every poset  $(P, \leq)$  is a simplicial complex $\Delta(P, \leq)$ called the {\em order complex} of the poset. $\Delta(P, \leq)$ has the set $P$ as vertices, and the finite chains of $(P, \leq)$ as faces. The poset homology of $(P,\leq)$ is by definition the simplicial homology of its order complex $\Delta(P, \leq)$.

\subsection{Taylor's resolution}
Minimal free resolutions of monomial ideals are very hard to compute in general. To solve this problem, two main lines of research have been adopted. One is to construct minimal resolutions for particular families of ideals, such as stable ideals \cite{EK90} or Borel-plus-powers \cite{M08} for example. In the context of system reliability, the resolutions of ideals of particular systems, such as $k$-out-of-$n$ systems and variants were studied in \cite{PS21,PSW21}. Another line of research consists on the construction of non-minimal resolutions that are easy to build, and which can then be minimalized following some standard procedure (like for instance the one described in \cite{CLO05}, Ch. 6, Theorem 3.15). Taylor's resolution \cite{T66} was the first such resolution, and the main one. Other general families of non-minimal resolutions include cellular resolutions \cite{BS98} or mapping cones \cite{CE95}. Taylor's resolution is strongly related to the $\lcm$-structure of the corresponding monomial ideal.
\begin{definition}
Let $I\subseteq R=\kb[x_1,\dots,x_n]$ be a monomial ideal minimally generated by $G(I)=\{\xb^{\mu_1},\dots,\xb^{\mu_r}\}$. The following sequence of modules and morphisms is a free resolution of $I$, and it is called the Taylor resolution of $I$:
\begin{equation}\label{eq:TaylorResolution}
\TT(I): 0\longrightarrow R^{{r}\choose{r}}\stackrel{\partial_r}{\longrightarrow}R^{{r}\choose{r-1}}\stackrel{\partial_{r-1}}{\longrightarrow}\cdots \longrightarrow R^{{r}\choose{2}}\stackrel{\partial_1}{\longrightarrow}R^{{r}\choose{1}}\longrightarrow I\longrightarrow 0,
\end{equation}
where the generators $e^i_{1},\dots,e^i_{{r}\choose{i+1}}$ of the $i$-th module in $\TT(I)$ are indexed by all the $\sigma\subseteq[r]$ such that $\vert\sigma\vert=i+1$. The multi-degree of each of these generators is $\mu_\sigma$. The differentials $\partial_i$ have the form
\[
\partial(e^i_\sigma)=\sum_s (-1)^{\sign(s,\sigma)} \frac{\xb^{\mu_\sigma}}{\xb^{\mu_{\sigma\setminus s}}} e^{i-1}_{\sigma\setminus s},
\]
where ${\sign(s,\sigma)}$ is the position of index $s$ in $\sigma$.
\end{definition}

Observe that since $\TT$ has a generator for each subset of $[r]$, the corresponding formula (\ref{eq:resHilbert}) for the $j$-reliability of a system $\Sc$ based on the form of $\HN(I_j(\Sc))$ obtained from $\TT(I_j(\Sc))$ corresponds to the usual inclusion-exclusion formula. Also, the structure (\ref{eq:TaylorResolution}) of the Taylor resolution is isomorphic to the chain complex associated to a simplicial complex, precisely the full simplex of dimension $r-1$. This fact is expressed by saying that the Taylor resolution is supported on the $(r-1)$-simplex. In general, monomial resolutions supported on simplicial and cellular complexes are called {\em cellular resolutions} \cite{BS98}.

\subsection{Hochster's formulas}\label{sec:Hochster}
Both monomial ideals and simplicial complexes have a notion of Betti numbers related to them. These notions are strongly connected through variations of the celebrated Hochster's formula.
\begin{theorem}[\cite{H77}]
Let $\Delta$ be a simplicial complex on $n$ vertices, and let $I_\Delta$ its Stanley-Reisner ideal. For any subset $W\subseteq[n]$ let $\Delta_W=\{\sigma\in\Delta\vert\sigma\subseteq W\}$. Then,
\begin{equation}
    \beta_{i,\xb_W}(I_\Delta)=\dim_\kb \tilde{H}_{\vert W\vert-i-2}(\Delta_{W};\kb),
\end{equation}
where $\xb_W=\prod_{i\in W}x_i$, and $\kb$ is a field.
\end{theorem}
This formula establishes the direct relation between Betti numbers of monomial ideals and those of simplicial complexes via the Stanley-Reisner correspondence. Variations of this formula exist for the Koszul simplicial complex and ideal, and for the poset homology of the lcm-lattice of monomial ideals.

\begin{theorem}\label{th:HochsterKoszul}
Let $\Delta$ be a simplicial complex on $n$ vertices and let $I^\Delta\subseteq \kb[x_1,\dots,x_n]$ its Koszul ideal, then
\begin{equation}
\dim \tilde{H}_i (\Delta; \kb) = \beta_{i+1,\1}(I^{\Delta})\;\;, \forall i\geq -1,    
\end{equation}
where $\1=x_1\cdots x_n$ is the monomial corresponding to the product of all variables.
\end{theorem}

\begin{theorem}\label{th:LatticeHomology}
    Let $I\subseteq \kb[x_1,\dots,x_n]$ be a monomial ideal and $L_I$ its lcm-lattice, then, for $\mu\in L_I$ 
    \begin{equation}
    \beta_{i,\mu}(I)=\dim_\kb \tilde{H}_{i-1}(\Delta(\hat{0},\mu)_{L_I};\kb),
\end{equation}
where $\supp(\mu)=\{i\in[n]\mid x_i \mbox{ divides } \xb^\mu\}$.
\end{theorem}
\section{The domination function as Euler characteristic}\label{sec:dominationEuler}
The study of domination theory can benefit from the correspondences among combinatorial objects described in the previous sections. Since a coherent systems has a monomial ideal associated to it, we can use tools for the study and computation of invariants of monomial ideals (and hence simplicial complexes and lattices) to study coherent systems and their domination functions. A first result is that, given a coherent system $\Sc$, the maximum-closure of its set of minimal $j$-paths is isomorphic to the $\lcm$-lattice of the corresponding $j$-reliability ideal.
\begin{theorem}\label{th:lcm_closure}
 Let $\Sc$ be a multi-state coherent system and $\overline{\Fc}_j(\Sc)$ its set of minimal $j$-paths. Let $I_j(\Sc)$ be the $j$-reliability ideal of $\Sc$. Then 
 \[
 L_{I_j(\Sc)}\cong \cl(\overline{\Fc}_j(\Sc)).
 \]
\end{theorem}
\begin{proof}
The ideal $I_j(\Sc)$ is minimally generated by the monomials corresponding to the minimal $j$-paths of $\Sc$, i.e
\[
G(I_j(\Sc))=\left\lbrace \xb^\mu\vert\mu \in \overline{\Fc}_j(\Sc)\right\rbrace,
\]
hence there is a bijection between $G(I_j(\Sc))$ and $\overline{\Fc}_j(\Sc)$.

Let now $\xb^\mu\in L_{I_j(\Sc)}$, then there exist some generators $\xb^{\mu_1},\dots,\xb^{\mu_i}$ such that $\xb^\mu=\lcm(\xb^{\mu_1},\dots,\xb^{\mu_i})$ and by definition $\mu=\mu_1\lor\dots\lor\mu_i$, therefore the bijection between $G(I_j(\Sc))$ and $\overline{\Fc}_j(\Sc)$ extends to a bijection between $L_{I_j(\Sc)}$ and $\cl(\overline{\Fc}_j(\Sc))$. This bijection preserves the ordering, and hence $L_{I_j(\Sc)}$ and $\cl(\overline{\Fc}_j(\Sc))$ are isomorphic as posets.
\end{proof}

Based on this result, the domination function of a coherent system, can be expressed as the Euler characteristic of the corresponding monomial ideal.

\begin{theorem}\label{th:dominationBetti}
Let $\Sc=(\Cc,\phi)$ be a coherent system, $\overline{\Fc}_j(\Sc)$ its set of minimal $j$-paths for some $j$, $I_j(\Sc)$ the $j$-reliability ideal of $\Sc$, and $\delta_j$ the domination function of $\cl(\overline{\Fc}_j(\Sc))$. Then
\[
    \delta_j(\mu)=\sum_i (-1)^i\beta_{i,\mu}(I_j(\Sc))=\chi_\mu(I_j(\Sc)).
\]
\end{theorem}
\begin{proof}
Let $\TT:=\TT(I_j(\Sc))$ be the Taylor resolution of the $j$-reliability ideal of $\Sc$. For each module $T_i$ in $\TT$ we denote by $T_{i,\mu}$ the submodule of $T_i$ generated in multi-degree $\mu$, and by $\gamma^T_{i,\mu}(I_j(\Sc))=\rank(T_{i,\mu})$ the rank of this submodule.
Each generator of $T_{i,\mu}$ corresponds to a set $\sigma\subseteq\{\xb^{\mu_1},\dots,\xb^{\mu_r}\}$ of the minimal generators of $I_j(\Sc)$, of cardinality $i+1$, such that $\xb_\sigma=\xb^\mu$. By the isomorphism of Theorem \ref{th:lcm_closure}, this generator corresponds to a formation of $\mu$ of cardinality $i+1$, and hence the alternating sum of the ranks of the Taylor resolution of $I_j(\Sc)$ at multi-degree $\mu$ equals the domination function $\delta_j(\mu)$.

Now, by the rank-nullity theorem, or by uniqueness of the Hilbert series, we know that the alternating sum of the ranks of any resolution of an ideal at a given multi-degree is the same, hence independent of the resolution. In particular, if we consider the minimal free resolution, then the ranks of its modules are called the (multi-graded) Betti numbers of the ideal, and denoted by $\beta_{i,\mu}$, i.e, we have
\[
\sum_i (-1)^i\gamma^T_{i,\mu}(I_j(\Sc))=\sum_i (-1)^i\beta_{i,\mu}(I_j(\Sc)).
\]

The alternating sum of the Betti numbers is, by Definition \ref{def:EulerCharIdeals}, the Euler characteristic of the ideal at multi-degree $\mu$, $\chi_\mu(I_j(\Sc))$, and hence the result.
\end{proof}

\begin{corollary}
    Let $\nu_j=\lcm\{\xb^\mu\mid\xb^\mu \in G(I_j(\Sc)\}$, then $\overline{d}_j(\Sc)=\chi_{\nu_j}(I_j(\Sc))$.
\end{corollary}

\noindent{\bf Notation: } We will use the notation $\delta_j(\Sc)$ and $d_j(\Sc)$ to express the corresponding Euler characteristics, thereby remarking the connection to the signed domination function and the signed domination of the $j$-level structure function of system $\Sc$.

\begin{corollary}
    The signed domination function of a multidegree is the coefficient of the corresponding monomial in the Hilbert numerator of the $j$-reliability ideal, i.e.
    \[
    HN_{I_j(\Sc)}(\xb)=\sum_\mu \delta_j(\mu)\xb^\mu.
    \]
\end{corollary}

\section{Computing $\chi_\mu(I_j(\Sc))$}\label{sec:computationMVT}
Theorem \ref{th:dominationBetti} provides a way to compute the domination function of any coherent system $\Sc$ via the ranks of the modules in any free resolution of the ideal $I_j(\Sc)$. Computing the multi-graded ranks of free resolutions is in general a demanding task, in particular if we consider the {\em minimal} free resolution. Several algorithms have been proposed in the literature, since this is one of the main problems in combinatorial commutative algebra. One such fast algorithm is the one based on Mayer-Vietoris Trees \cite{S09}, which is, in the worst case, exponential on the number of generators of the ideal, i.e. on the number of minimal $j$-paths of the system. If the number of minimal paths of the system is smaller than the number of components, this approach is very efficient. However, if the number of minimal paths is larger than the number of components, then we can propose a more efficient approach for the computation of multi-graded Euler characteristic, which is exponential on the number of components of the system (i.e. variables of the ring). These two approaches can furthermore provide insight on the structure of the domination and domination function hence give some characterizations and shortcuts for certain types of systems, see Section \ref{sec:examples}.
\subsection{Computation of the domination function using Mayer-Vietoris trees}
Among the general approaches to the computation of free resolutions of monomial ideal, a very important one uses cone resolutions \cite{CE95}. The support of these resolutions can be computed in an efficient way without needing to compute the differentials using Mayer-Vietoris trees \cite{S09}. This algorithm was adapted in \cite{BHS19} to compute simplicial homology by the computation of the $\1$-strand of a cone resolution (i.e. the Betti numbers $\beta_{i,\1}(I)$ at multi-degree $\1=(1,\dots,1)$). Here we introduce the basics of Mayer-Vietoris trees and we adapt them to compute (bounds of) $\beta_{i,\mu}(I)$ for any $\mu$ and all $i\in\{0,\dots,n\}$.

Let $I\subseteq \kb[x_1,\dots,x_n]$ be a monomial ideal and $G(I) = \{\xb^{\mu_1}, . . . ,\xb^{\mu_r}\}$ a (non necessarily minimal) monomial generating set of $I$. Let $I_i=\langle \xb^{\mu_1}, . . . ,\xb^{\mu_i}\rangle$ be the ideal generated by the first $i$ generators of $I$. For each $i$ consider the following exact sequence
\[
0\rightarrow I_{i-1}\cap\langle\xb^\mu_i\rangle\stackrel{\jk}{\rightarrow} I_{i-1}\oplus\langle\xb^\mu_i\rangle\rightarrow I_i\rightarrow 0.
\]
Assuming that we know free resolutions $\FF'_i$ and $\tilde{\FF}_i$ for $I'_i=I_{i-1}$ and $\tilde{I}_i=I_{i-1}\cap\langle \xb^{\mu_i}\rangle$ respectively, then a (not necessarily minimal) resolution $\FF_i$ of $I_i$ is given by the so called mapping cone of the chain complex morphism $\psi: \tilde{\FF}_i\rightarrow \FF'_i$ that lifts the inclusion $\jk$. Using this sequence recursively on $i$ we can compute a free resolution of $I$ that is called an {\em iterated mapping cone resolution}. The ideals involved in this iterative process can be displayed as a binary tree in which each node is given by an ideal. The root of the tree is $I$ and every node $J = \langle \xb^{\nu_1},\dots, \xb^{\nu_j}\rangle$ has the ideal $J'=\langle \xb^{\nu_1},\dots, \xb^{\nu_{j-1}}\rangle$ as right child, and $\tilde{J}=J'\cap\langle\xb^{\nu_j}\rangle$ as left child. This is called a {\em Mayer-Vietoris tree} of $I$, denoted $MVT(I)$, see \cite{S09} for details.

Each node in a Mayer-Vietoris tree is assigned a position and a dimension. The root has position $1$ and dimension $0$, the right child of the node at position $p$ and dimension $d$ is has position $2p+1$ and dimension $d$, and its left child has position $2p$ and dimension $d+1$. We denote by $MVT(I)_{p,d}$ the set of monomials generating the ideal in the node of $MVT(I)$ at position $p$ and dimension $d$. A node is called {\em relevant} if it is either the root of the tree, or if it is at an even position. Using this notation, the following holds.
\begin{proposition}\label{prop:BettiMVT}
 Let $I\subseteq \kb[x_1\dots,x_n]$ be a monomial ideal and $\FF(I)$ an iterated cone resolution of $I$ whose support is given by the Mayer-Vietoris tree $MVT(I)$. Let $F_i$ denote the $i$-th module of $\FF$. Then
 \[
 \rank_\mu(F_i)=\#(\{p \mbox{ even, or }p=1\mid \xb^\mu\in MVT(I)_{p,i}\}).
 \]
\end{proposition}
\begin{proof}
    The proof is simply the fact that the Mayer-Vietoris tree describes the support of the corresponding resolution. This support is given by the relevant nodes in the tree.
\end{proof}

The computation of $\chi_\mu(I)$ using Proposition \ref{prop:BettiMVT} amounts then to obtain all the nodes in which $\xb^\mu$ is present and take the alternating sum according to the dimension of the nodes. The usual construction of the tree can be optimized so that the branches in which $\xb^\mu$ is not present, are not considered. This can be done efficiently using adequate orderings of the generators of each node (equivalently, choosing good pivots at each node to construct the tree). In the case $\mu=\1$ we select as pivots those monomials having the least frequent variable, so that we obtain nodes that do not have full support, and hence $\xb^\1$ is not present in the corresponding branch. This is the strategy followed in \cite{BHS19} for the fast computation of Betti numbers of simplicial complexes in high dimension, see Example \ref{ex:MVTEuler} (a). In the general case, another good strategy is to chose as pivots those monomials that divide $\xb^\mu$ so that we obtain as soon as possible those nodes that do not divide $\xb^\mu$ and hence the corresponding branches can be pruned, see Example \ref{ex:MVTEuler} (b).

The size of $MVT(I)$ is, in the worst case $2^r-1$, where $r$ is the number of generators of $I\subseteq\kb[x_1,\dots,x_n]$, and usually much smaller. Hence, in the cases when $r<n$ this is a convenient approach to an efficient computation of $\delta_j(\mu)=\chi_\mu(I_j(\Sc))$. 

\begin{example}\label{ex:MVTEuler} {\rm (a)} Consider a multi-state system $\Sc=(C,\phi)$ with $n=7$ components. The state spaces of the components are $S_1=S_2=\{0,1,2,3\}$ and $S_3=\cdots=S_7=\{0,1\}$. The structure function of the system is $\phi(\cb)=\sum_{i=1}^7c_i$, hence $S=\{0,\dots,11\}$. The ideal $I=I_{10}(\Sc)$ is minimally generated by
 \begin{align*}
 I=&\{x_1^3x_2^3x_3x_4x_5x_6,x_1^3x_2^3x_3x_4x_5x_7,x_1^3x_2^3x_3x_4x_6x_7,x_1^3x_2^3x_3x_5x_6x_7,\\
 &x_1^3x_2^3x_4x_5x_6x_7,x_1^3x_2^2x_3x_4x_5x_6x_7,x_1^2x_2^3x_3x_4x_5x_6x_7\}.
 \end{align*}

Taking as pivot the last generator of each node, the relevant nodes of $MVT(I)$ are:
\begin{align*}
&MVT_{2,1}(I)=MVT_{6,1}(I)=MVT_{14,1}(I)=MVT_{30,1}(I)=\\
&MVT_{62,1}(I)=MVT_{126,1}(I)=\langle x_1^3x_2^3x_3x_4x_5x_6x_7\rangle,\\
\end{align*}
and $MVT_{1,0}(I)=G(I)$. Hence $\delta_{10}(\xb)=1$ for all minimal $10$-paths of the system, $\xb\in\overline{\Fc}_{10}(\Sc)$. The signed domination of $\phi_{10}$ for the path $(3,3,1,1,1,1,1)$ is $\delta_{10}(3,3,1,1,1,1,1)=-6$, and $\delta_{10}(\xb)=0$ for every other state vector of the components of the system.

{\rm (b)} Consider a binary system $\Sc=(C,\phi)$ with $5$ components, such that the set of minimal paths is
 
\begin{align*}
\overline{\Fc}(\Sc)=&\{(1,0,1,1,0),(1,1,1,0,0),(0,1,0,1,0),(0,0,1,0,1),\\
&(1,0,0,0,1),(0,0,0,1,1)\}.
\end{align*}

The corresponding ideal is
\[
I(\Sc)=\langle xzt,xyz,yt,zu,xu,tu\rangle\subseteq\kb[x,y,z,t,u].
\]
In order to compute $\delta({1,0,0,1,1})$, which corresponds to the monomial $xtu$, we use as pivots the generators of $I$ that divide $xtu$, namely $tu$ and $xu$. We obtain the following tree (using as pivot at each step the last monomial in the node), pruned at the nodes that have no divisors of $xtu$ .
\begin{center}
 \begin{tikzpicture}[scale=1,transform shape]
 \tikzstyle{level 1}=[sibling distance=4.5cm]
 \tikzstyle{level 2}=[sibling distance=3.5cm]
 \node{$(1,0)$ xzt,xyz,yt,zu,xu,tu}
 child{ node{$(2,1)$ $ytu,ztu,\mathbf {xtu}$}}
 child{ node[color=black!50!white]{$(3,0)$ $xzt,xyz,yt,zu,xu$}
 	child{node{$(6,1)$ $xzu,xytu$}}
 	child{node[color=black!50!white]{$(7,0)$ $xzt,xyz,yt,zu$}
	}};
 \end{tikzpicture}
 \end{center}

We see from the tree that $\delta(1,0,0,1,1)=-1$.
\end{example}
\subsection{Computation of the domination function using Koszul complexes}\label{sec:computationKoszul}
Recall that for any monomial ideal $I\subseteq\kb[x_1,\dots,x_n]$ and multi-degree $\mu\in\NN^n$ we can define a simplicial complex on $\vert \supp(\mu)\vert\leq n$ vertices, $\Delta K(I,\mu)$, called the {\em Koszul simplicial complex} of $I$ at $\mu$, see Definition \ref{def:KoszulComplex}. By Hochster's formula, Theorem \ref{th:HochsterKoszul}, applied to $I=I_j(\Sc))$, we have that
\[
\chi_\mu(I_j(\Sc))=\sum_i (-1)^i\beta_{i,\mu}(I_j(\Sc))=\sum_i (-1)^i\dim_\kb \tilde H_{i-1}(\Delta K(I_j(\Sc),\mu);\kb),
\]
hence the domination and domination function of $\Sc$ equals minus the reduced Euler characteristic of the Koszul complex of $I_j(\Sc)$.
\begin{theorem}\label{th:dominationKoszul}
Let $\Sc$ be a system, $\delta_j$ its signed domination function for level $j$, and $\overline d_j(\Sc)$ its level domination for level $j$. Then
\[
\delta_j(\mu)=-\tilde\chi(\Delta K(I_j(\Sc),\mu)) ,\mbox{ for } \mu\in\NN,
\]
and
\[
\overline{d}_j(\Sc)=-\tilde\chi(\Delta K(I_j(\Sc),\nu_j)), 
\]
where $\nu_j=\lcm\{\mu\mid \mu \in \overline{\Fc}_j(\Sc)\}$.
\end{theorem}

\begin{remark}
Observe that
\[
  d_j(\Sc)=\begin{cases}\overline d_j(\Sc)&\text{if }\nu_j=\mb,\\ 0&\text{otherwise.}\end{cases}
\]    
\end{remark}

\begin{example}
 Consider the system in Example 4.6 in \cite{H25}. The system has $4$ components, each of which can be in any of the states $\{0,1,2\}$, and the structure function of the system is given by
 \[
 \phi(\xb)=\sum_{i=1}^4 x_i,
 \]
hence the system itself can be in any of the levels $\{0,1,\dots,8\}$. We consider the $4$-level structure function
\[
\phi_4(\xb)=I\left( \sum_{i=1}^4x_i\geq 4\right).
\]
The ideal $I_4(\Sc)$ is minimally generated by $19$ monomials, hence the Taylor resolution has $2^{19}$ generators in total (which is also the worst case size of a Mayer-Vietoris tree of $I_4(\Sc)$). In order to compute the domination of $\phi_4$ we only need to check the multi-graded Euler characteristic at multi-degree $(2,2,2,2)$. Since $\xb^{(1,1,1,1)}\in I_4(\Sc)$, the Koszul simplicial complex $\Delta K_{(2,2,2,2)}(I_4(\Sc))$ is the full simplex on $4$ vertices, hence its reduced Euler characteristic is $0$ and therefore $d_4(\Sc)=0$.
\end{example}

This example can be generalized to the following Corollary to Theorem \ref{th:dominationKoszul}.
\begin{corollary}\label{cor:criterionM}
Let $\phi_j$ be the structure function of a MMS system at level $j$, and let the sets of the states of the components be given by $\Sc_i=\{0,1,\dots,m_i\}$. If $\phi_j(m_1-1,\dots,m_n-1)=1$, then $d_j(\Sc)=0$.
\end{corollary}
\begin{proof}
    The domination $d_j(\Sc)$ is given by the reduced Euler characteristic of the simplicial complex $\Delta K_{m_1,\dots,m_n}(I_j(\Sc))$. Since $\phi_j(m_1-1,\dots,m_n-1)=1$, then we have that $\xb^{(m_1,\dots,m_n)}/\xb_\sigma\in I_j(\Sc)$ for every subset $\sigma$ of $[n]$. Hence $\Delta K_{m_1,\dots,m_n}(I_j(\Sc))$ is the full simplex in $n$ variables, whose reduced Euler characteristic is $0$.
\end{proof}

\begin{corollary}
Let $\phi_j$ be the structure function of a MMS system at level $j$, and $\nu_j=\lcm\{\mu\mid \mu \in \overline{\Fc}_j(\Sc)\}$. If $\phi_j({\nu_j}_1-1,\dots,{\nu_j}_n-1)=1$, then $\overline d_j(\Sc)=0$.
\end{corollary}
\begin{proof}
    The proof follows the same arguments as in Corollary \ref{cor:criterionM}
\end{proof}
\begin{remark}
Observe that $\Delta K(I_{\Sc},\1)$ is the simplicial complex whose facets are the complements of the sets corresponding to the minimal paths of the system, and every simplicial complex arises in this way. Hence computing $d(\phi)$ from the minimal paths is equivalent to the Euler characteristic problem for complexes given by facets, which is \#P-complete, as shown in \cite{RS13}.
\end{remark}

\subsection{A note on the M\"obius inversion formula}
The classical M\"obius inversion formula $d(A)=\sum_{B\subseteq A}(-1)^{|A\setminus B|}\phi(B)$ \cite{R64} is used in \cite{H25} to compute signed domination of binary and multi-state systems (see also \cite{B82}). The previous section contributes a topological interpretation of the relation between domination theory and M\"obius functions, that will be used in the next sections, together with the tools of algebraic reliability.
 
\begin{proposition}\label{prop:MobiusInversionKoszul}
For every monotone $\phi$, level $j$ and $\mu\in\mathbb \Cc$,
\[
  \delta_j(\mu)=\sum_{\sigma\subseteq\supp\mu}(-1)^{|\sigma|}\phi_j(\mu-{\bf 1}_\sigma)
  =\sum_{\sigma\in\Delta K(I_j,\mu)}(-1)^{|\sigma|}=-\tilde\chi\bigl(\Delta K(I_j,\mu)\bigr).
\]
\end{proposition}
 
\begin{proof}
The first equality is M\"obius inversion on the product of chains, whose
M\"obius function is $\prod_i{\mathbf{M}}(a_i,b_i)$ with $\mathbf{M}(a,a)=1$,
$\mathbf{M}(a,a+1)=-1$ and $0$ otherwise. The second holds because
$\phi_j(\mu-{\bf 1}_\sigma)=1\iff x^{\mu}/x^\sigma\in I_j\iff\sigma\in\Delta K(I_j,\mu)$.
The third is the definition of $\tilde \chi$.
\end{proof}

Furthermore, the relation between the homology of $\Delta K(I)$ and that of the $\lcm$-lattice $L_I$ (see Theorem \ref{th:LatticeHomology}) provides a formulation of the relation between domination and Euler characteristic in terms of the M\"obius function, via Rota's crosscut theorem \cite{R64} by seeing the signed domination function as minus the M\"obius function of the lattice of formations (see \cite{H25}).

\section{Multi-state systems via binary systems and vice-versa}\label{sec:polarization}
Given a general monomial ideal, the operation of {\em polarization} gives a squarefree monomial ideal that shares, among other invariants, the Betti numbers of the original one. This implies that in the particular case of reliability ideals, we can polarize the $j$-reliability ideal of a multi-state coherent system and obtain a coherent ideal which corresponds to a certain binary system, that can then be called the polarization of the $j$-level structure function of the original system. Therefore we have a way to study multi-state systems via binary systems. In general, the number of components of the binary system is bigger than the number of components of the original system. This is an algebraic counterpart of the proposed decomposition theorem for multi-state structure functions in \cite{BS82}, see also \cite{W85}.

Conversely, given a squarefree ideal, the operation of {\em depolarization} gives a family of non-squarefree ideals whose polarizations are isomorphic to the original ideal, and hence share its Betti numbers. These depolarizations are ideals in a smaller number of variables than the original squarefree one. We therefore have a way to study certain properties of binary systems using multi-state systems in a smaller number of components. In particular, this is a convenient way to study the signed domination function of binary systems. The main advantage is to reduce the number of variables, and hence one can make use of the approach in Section \ref{sec:computationKoszul} for instance.

We give here a brief account of these two operations on monomial ideals and examples of application to signed domination theory. For further details, see \cite{F05,HH11,MMVV19,MPSW20}.

\subsection{Polarization}
 \begin{definition}\label{def:polarization}
Let us consider two multi-indices $\ab=(a_1,\dots,a_n)$ and $\mu=(\mu_1,\dots,\mu_n)$ in $\NN^n$, such that $\mu_i\leq a_i$ for all $i$. The {\em polarization} of $\mu$ in $\NN^{a_1+\cdots+a_n}$ is the multi-index 
\[
\Pc(\mu)=(\underbrace{1,\dots,1}_{\mu_1},\underbrace{0,\dots,0}_{a_1-\mu_1},\dots,\underbrace{1,\dots,1}_{\mu_n},\underbrace{0,\dots,0}_{a_n-\mu_n}).
\]
The {\em polarization} of the monomial $\xb^\mu=x_1^{\mu_1}\cdots x_n^{\mu_n}\in R=\kb[x_1,\dots,x_n],$ with respect to $\ab$ is the squarefree monomial 
\[
\Pc(\xb^\mu)=\xb^{\Pc(\mu)}=x_{1,1}\cdots x_{1,\mu_1}\cdots x_{n,1}\cdots x_{n,\mu_n}
\]
in the ring $\Pc(R)=\kb[x_{1,1},\dots,x_{1,a_1},\dots,x_{n,1},\dots,x_{n,a_n}]$.

Let now $I=\langle \xb^{\mu_1},\dots,\xb^{\mu_r}\rangle\subseteq R$ be a monomial ideal and let $\xb^{\mu_I}=\lcm(\xb^{\mu_1},\dots,\xb^{\mu_r})$.
The {\em polarization of $I$}, $\Pc(I)$, is the monomial ideal in $\Pc(R)$ given by $\Pc(I)=\langle \Pc(\xb^{\mu_1}),\dots,\Pc(\xb^{\mu_r})\rangle$, where $\Pc(\xb^{\mu_i})$ is the polarization of $\xb^{\mu_i}$ with respect to $\mu_I$, for each $i$.
\end{definition}

\begin{example}
    Let $I=\langle x^3y,y^2z,z^2\rangle\subseteq\kb[x,y,z]$ a monomial ideal. Then 
    \[
    \Pc(I)=\langle x_1x_2x_3y_1,y_1y_2z_1,z_1z_2\rangle\subseteq\kb[x_1,x_2,x_3,y_1,y_2,z_1,z_2].
    \]
\end{example}

The polarization operation has been widely studied, mainly because of the features that a monomial ideal $I$ and its polarization $\Pc(I)$ share, see e.g. \cite{F05,HH11,V15,AV22}. 
Among the important properties shared by a monomial ideal and its polarization, we are mainly interested in Betti numbers (see also \cite{F05,MMVV19,AFL22}).

\begin{proposition}[see Corollary 1.6.3, item (1) in \cite{HH11}]\label{prop:BettiPolarization}
Let $I\subseteq R=\kb[x_1,\dots,x_n]$ be a monomial ideal and let $\Pc(I)\subseteq \Pc(R)$ be its polarization. Then
\[
\beta_{i,\mu}(I)=\beta_{i,\Pc(\mu)}(\Pc(I)) \mbox{ for all }i\in\ZZ \mbox{ and }\mu\in\NN^n.
\]
\end{proposition}

\subsection{Depolarization}\label{sec:depolarization}
The {\em depolarization} of a squarefree monomial ideal $J$ is not unique, in the sense that it gives a family of monomial ideals with the main homological features of $J$, but in rings with smaller number of variables, see \cite{MPSW20,PS22}. When applied to simplicial complexes, this allows to study them via other complexes in a smaller number of vertices.

\begin{definition}
Let $R,S$ and $T$ be polynomial rings over a field $\kb$. Let $I\subseteq R$ be a squarefree monomial ideal. A {\em depolarization of $I$} is a monomial ideal $J\subseteq S$ such that $I$ is isomorphic to $\Pc(J)\subseteq T\cong R$, i.e., there is a bijective map $\varphi$ from the variables of $R$ to the variables of $T$ such that $\varphi(G(I))=G(\Pc(J))$.
\end{definition}



The set of depolarizations of a given squarefree monomial ideal has an order structure, that can be fully described by the chains of the support poset of the ideal, see \cite{PS22} for details.

\begin{example}
 Consider a binary system $\Sc$ with $6$ components that works whenever at least $3$ consecutive components are working. Its squarefree monomial ideal is 
 \[
 I({\Sc})=I=\langle x_1x_2x_3,x_2x_3x_4,x_3x_4x_5,x_4x_5x_6\rangle\subseteq\kb[x_1,\dots,x_6].
 \]
A depolarization of $I$ is the ideal
\[
J=\langle a^3,a^2b,ab^2,b^3\rangle \subseteq\kb[a,b].
\]
From Proposition \ref{prop:BettiPolarization} we have that $\chi_{\Pc(\mu)}(I)=\chi_\mu(J)$ for all $\mu$, hence we can use $J$ to compute the signed domination of $\Sc$. Since $J$ has only $2$ variables, using the results in Section \ref{sec:computationKoszul} to compute $\delta(\mu)$ for each $\mu$ needs a simplicial complex in just $2$ vertices. We have that
{\small
\[   
\delta({\mu})= 
     \begin{cases}
   \;\;\,1 &\quad\text{if }\mu\in\{(1,1,1,0,0,0),(0,1,1,1,0,0),(0,0,1,1,1,0),(0,0,0,1,1,1)\}\\
       -1&\quad\text{if }\mu\in\{(1,1,1,1,0,0),(0,1,1,1,1,0),(0,0,1,1,1,1)\}\\
       \;\;\,0&\quad\text{otherwise.} \\ 
     \end{cases}
\]
}
\end{example}

\section{Signed domination of certain families of coherent systems}\label{sec:examples}
\subsection{Dominant systems and systems with dominant components}
A notion of dominant system is given in \cite{HZ04}. Under this notion, a dominant system is such that all minimal path vectors to state $j$ cause the system to be exactly in state $j$. I.e. $\phi$ is \emph{dominant} if
$\phi(\yb)>\phi(\xb)$ implies that either $\yb>\xb$ or $\yb>\zb$ for
some $\zb$ with $\phi(\zb)=\phi(\xb)$. Dominance, in this setting, has a characterization in terms of the supports
of the level-wise signed domination functions, simply observing that the family of minimal path sets of the levels is pairwise disjoint, and that first-order cancellations do not occur at any level. It constrains only the relationship \emph{between} levels and does not simplify the computation of signed domination functions \emph{within} a level. We therefore propose a different notion of dominance, closer to the one in \cite{A17}, which works level-wise and is based on how a certain individual component performance level dominates or not the performance level of other components in the minimal $j$-paths of the system for a given $j$. Also, we give in Proposition \ref{prop:algebraic RelDominant} below an algebraic characterization of dominant systems in the traditional sense of \cite{HZ04}.

Let $\Sc=(C,\phi)$ be a multi-state component system such that the state set of each component $\overline{c}_i$ is given by $S_i=\{0,\dots,m_i\}$.
\begin{definition}\label{def:dominant}
We say that a minimal $j$-path $\xb \in \overline{\Fc}_j$ is {\em dominant} for component $i$ if $c_i=m_i$ in $\xb$ and $c_i\lneq m_i$ in every other $\yb \in \overline{\Fc}_j$.

We say that $\Sc=(C,\phi)$ is a {\em component-based dominant} system for level $j$ if every minimal $j$-path is dominant for at least one component.
\end{definition}
 In the rest of the section, we refer to component-based dominant systems, simply by dominant.
    Observe that by definition, a dominant system has at most one dominant minimal $j$-path for each component, and a minimal $j$-path can be dominant for several components.

\begin{example}\label{ex:dominantSystem}
    Let $\Sc=(C,\phi)$ be a multi-state system with $4$ components, such that
    \[
        S_1=S_3=S_4=\{0,1,2\}\text{ and }S_2=\{0,1,2,3\}.
    \]
    Let $\phi_j$ for some $j$ be given by the following minimal $j$-paths
    \[
    \overline{\Fc}_j=\{\xb_1,\xb_2,\xb_3\}=\{(2,1,2,1),(0,2,0,2),(1,3,0,0)\}.
    \]
    Then, $\Sc$ is a dominant system at level $j$, $\xb_1$ is dominant for components $1$ and $3$, $\xb_2$ is dominant for component $4$, and $\xb_3$ is dominant for component $2$. 
\end{example}

Using the results in Section \ref{sec:computationMVT}, we can compute the signed domination function for a dominant system using Mayer-Vietoris trees.

\begin{proposition}
Let $\Sc=(C,\phi)$ be a dominant system for level $j$ and let $\overline{\Fc}_j=\{\xb_1,\dots,\xb_r\}$. Then 
\[   
\delta_j({\mu})= 
     \begin{cases}
       (-1)^{k-1}&\quad\text{if }\xb_\mu=\xb_{i_1}\vee\cdots\vee \xb_{i_k} \mbox{ for a subset } \{\xb_{i_1},\dots,\xb_{i_k}\}\subseteq\overline{\Fc}_j\\
       \;\;\,0&\quad\text{otherwise.} \\ 
     \end{cases}
\]
In particular, $\overline d_j(\Sc)=(-1)^{r-1}$.
\end{proposition}
\begin{proof}
 Consider a Mayer-Vietoris tree of the ideal $I_j(\Sc)=\langle \xb_1,\dots,\xb_r\rangle$ in which each generator is the monomial correspondent to each of the minimal $j$-paths of $\Sc$. Take any of them as pivot (for ease of notation, we assume that at each step the last monomial is the pivot).
 
 First, observe that if $r=1$ the result is obvious. We will proceed by induction on $r$, assuming that the result holds for every dominant system with $r-1$ minimal paths. We have that $MVT_{1,0}(I_j(\Sc))=I_j(\Sc)$. Now, take $\xb_r$ as pivot and let $\{i^1_r,\dots,i^k_r\}$ the set of components (i.e. variables) for which $\xb_r$ is dominant. We obtain that $MVT_{3,0}(I_j(\Sc))=\langle \xb_1,\dots,\xb_{r-1}\rangle$, which is the ideal of a dominant system in $r-1$ minimal paths and such that the maximum state of the components $\overline{c}_{i^1_r},\dots,\overline{c}_{i^k_r}$ is strictly lower than $m_{i^1_r},\dots,m_{i^k_r}$ respectively, and hence there is no repeated multi-degree between the tree hanging from node $(3,0)$ and the tree hanging from node $(2,1)$ . On the other hand, $MVT_{2,1}(I_j(\Sc))=\langle \lcm(\xb_1,\xb_r),\dots,\lcm(\xb_{r-1},\xb_r)\rangle$, hence all the generators of $MVT_{2,1}(I_j(\Sc))$ have a common factor, namely $x_{i^1_r}^{m^1_r}\cdots x_{i^k_r}^{m^k_r}$. Now, taking away the common factor, we have that the ideal in this node has exactly $r-1$ generators in at most $n-k$ variables. This ideal corresponds to a dominant system, since we only removed for the common factor the variables for which $\xb_r$ is dominant. Observe that any subset $T$ having the generator $\xb_r$ corresponds to the subset $T\setminus\xb_r$ of the set of $r-1$ generators of the left child. Since the left child increases the dimension by one, we have that the contribution of this subset changes the sign of the contribution of the given subset. Hence the result holds by the induction hypothesis.
\end{proof}

\begin{example}\label{ex:dominantMVT}
The Mayer-Vietoris tree of the ideal of the system in Example \ref{ex:dominantSystem} is
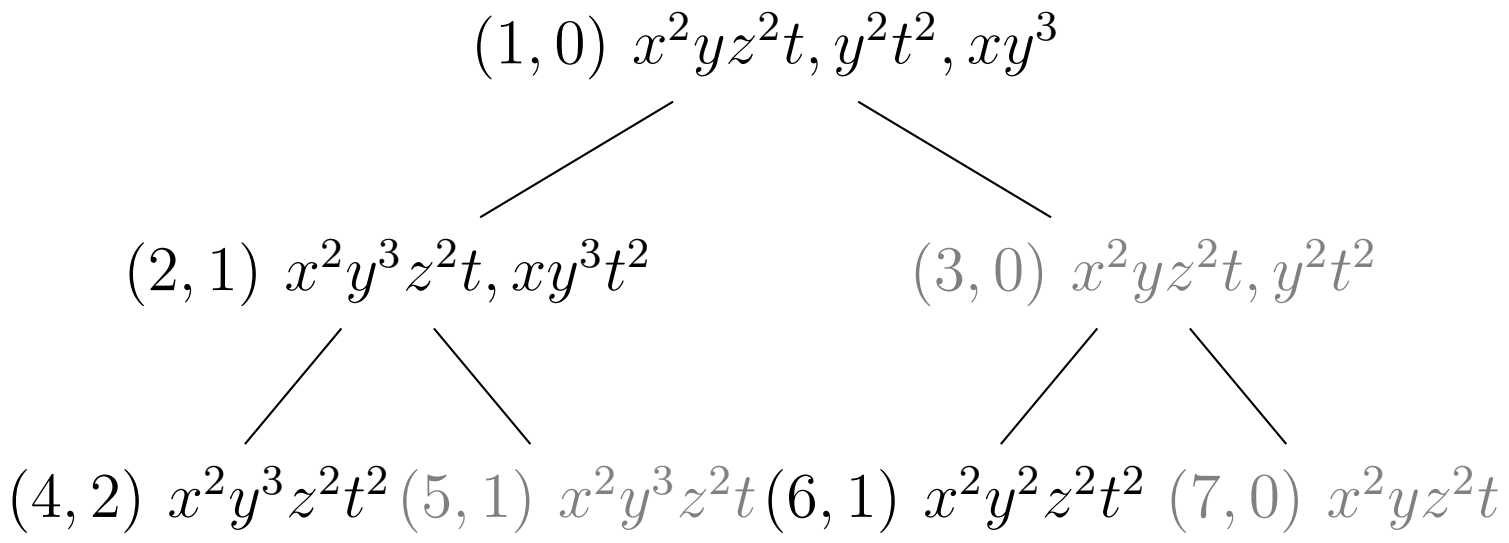
\begin{center}
 \begin{tikzpicture}[scale=1,transform shape]
 \tikzstyle{level 1}=[sibling distance=5cm]
 \tikzstyle{level 2}=[sibling distance=2.5cm]
 \node{$(1,0)$ $x^2yz^2t,y^2t^2,xy^3$}
 child{ node{$(2,1)$ $x^2y^3z^2t,xy^3t^2$}
 child{node{$(4,2)$ $x^2y^3z^2t^2$}}
 	child{node[color=black!50!white]{$(5,1)$ $x^2y^3z^2t$}
	}}
 child{ node[color=black!50!white]{$(3,0)$ $x^2yz^2t,y^2t^2$}
 	child{node{$(6,1)$ $x^2y^2z^2t^2$}}
 	child{node[color=black!50!white]{$(7,0)$ $x^2yz^2t$}
	}};
 \end{tikzpicture}
 \end{center}
From which
\begin{align*}
&\delta_j(2,1,2,1)=\delta_j(0,2,0,2)=\delta_j(1,3,0,0)=1,\\
&\delta_j(2,3,2,1)=\delta_j(1,3,0,2)=\delta_j(2,2,2,2)=-1,\\
&\delta_j(2,3,2,2)=d_j(\Sc)=1,\\
&\delta_j(\xb)=0\mbox{ otherwise.}
\end{align*}
\end{example}

\begin{proposition}\label{prop:algebraic RelDominant}
 A system $(\Sc,\phi)$ is dominant in the sense of \cite{HZ04} if and only if $I_{j+1}(\Sc)\subseteq\mathfrak{m}I_j(\Sc)$ for all $j$, where $\mathfrak{m}=\langle x_1,\dots,x_n\rangle$.
\end{proposition}
\begin{proof}
 First recall that $\phi(\xb)=j$ for all $\xb \in \overline\Fc_j(\Sc)$. In fact, let $\xb\in\overline\Fc_j(\Sc)$ such that $\phi(\xb)>j$. Since $(\Sc,\phi)$ is dominant, there exists some $\zb$ such that $\zb<\xb$ and $\phi(\zb)=j$, which contradicts the minimality of $\xb$, hence $\phi(\xb)=j$ for all $\xb\in\overline\Fc_j(\Sc)$. 

 
The fact that $\phi(\xb)=j$ for all $\xb\in\overline\Fc_j(\Sc)$ means that $G(I_j(\Sc))\cap I_{j+1}(\Sc)=\emptyset$. Since the monomials of $I_j(\Sc)$ outside $G(I_j)$ are those of $\mathfrak{m} I_j(\Sc)$, this means that $I_{j+1}(\Sc)\subseteq\mathfrak{m} I_j(\Sc)$. Conversely, the monomials in $G(\mathfrak{m} I_j(\Sc))$ have degree $j+1$, hence $\phi(\xb)=j+1$ for all $\xb\in\overline \Fc_{j+1}(\Sc)$.
\end{proof}

\subsection{Binary $k$-out-of-$n$ and consecutive linear $k$-out-of-$n$ systems.}\label{sec:bin_k-out-of-n}

\subsubsection{$k$-out-of-$n$ systems}\label{sec:k-out-of-n}
A \emph{binary $k$-out-of-$n$:G system} (G is for {\em good}), denoted by $\Sc_{(k,n)}$, is a binary system with $n$ components that works whenever at least $k$ components are working. Identically, we can denote by binary $k$-out-of-$n$:F system (F is for {\em fail}) a binary system with $n$ components that fails if at least $k$ components are not working. For more details about binary $k$-out-of-$n$ systems and their applications, we refer to \cite[Ch. 7]{KZ03}. We will work with $k$-out-of-$n$:G systems.

In \cite{SW09} the authors characterize binary $k$-out-of-$n$:G systems via monomial ideals. In fact, they provide explicit formulas for the Betti numbers of such ideals. As we are interested in computing the domination of these systems, we need to focus on the \emph{maximal element} \textbf{m}, which in the case of binary $k$-out-of-$n$ systems is $\textbf{m}=x_1\cdots x_n$.

The squarefree monomial ideal corresponding to a $k$-out-of-$n$ system is \[I_{(k,n)}=\left\langle \prod_{\substack{\sigma\subseteq\lbrace1,\dots,n\rbrace \\ \vert\sigma\vert=k}} x_i\ \vert\, i\in\sigma \right\rangle,\] i.e. all the possible products of any $k$ of the $n$ variables.

From Theorem \ref{th:dominationBetti} we have that 
\[
\delta(\textbf{m})=\sum_i (-1)^i\beta_{i,\textbf{m}}(I_{(k,n)}).
\]

The Betti numbers of $k$-out-of-$n$ systems were obtained in Section 5.1. of \cite{SW09}, where it is seen that if $\mu$ is any product of $k+i$ variables, then $\beta_{i,\mu}(I_{(k,n)})={{k+i-1}\choose{k-1}}$, and $\beta_{j,\nu}(I_{(k,n)})=0$ in all other cases. Hence we have the following (see \cite{AB84}).

\begin{proposition}\label{prop:dom_binaryKN}
Let $\Sc_{(k,n)}$ be a $k$-out-of-$n$ system. $\delta(\mu)=(-1)^i{{k+i-1}\choose{k-1}}$ for all products of $k+i$ variables, and  is $0$ for all other vectors. In particular $d(\Sc_{(k,n)})=(-1)^{n-k}{{n-1}\choose{k-1}}$.
\end{proposition}

\subsubsection{Binary consecutive linear $k$-out-of-$n$ systems.}
A \emph{binary consecutive linear $k$-out-of-$n$:G system}, denoted by $\Sc_{(k,n)}^{\text{con}}$, is a binary system with $n$ components that works if $k$ consecutive components are working \cite{KZ03}.

The squarefree monomial ideal corresponding to these systems is \[I_{(k,n)}^{\text{con}}=\langle x_1\dots x_k,x_2\dots x_{k+1},\dots,x_{n-k+1}\dots x_n \rangle,\] i.e. the ideal generated by any product of $k$ consecutive variables.
As in binary $k$-out-of-$n$ systems, we are interested in $\beta_{i,\1}(I_{(k,n)}^{\text{con}})$, where $\1=x_1x_2\cdots x_n$. We propose a recursive way to compute it:

\begin{proposition}\label{prop:BettiConsecBinKN}
	Let $\Sc_{(k,n)}^{\text{con}}$ be a binary consecutive linear $k$-out-of-$n$:$G$ system, with $0<k\leq n$, and $I_{(k,n)}^{\text{con}}$ its reliability ideal. Then, the domination of this system is given by
	 
		\centering
		$d(\Sc_{(k,n)}^{\text{con}}) = \begin{cases}
        1 & \text{if } n\equiv k\mod(k+1),\\
        -1  & \text{if } n\equiv 0\mod(k+1),\\
        0  & \text{otherwise}.
        \end{cases}$
\end{proposition}

\begin{proof}
Let $\Sc_{(k,n)}^{\text{con}}$ be a binary consecutive linear $k$-out-of-$n$:$G$ system, with $0<k\leq n$, and $I_{(k,n)}^{\text{con}}$ its reliability ideal. 

In order to prove this result we will show that the domination of this system, $d(\Sc_{(k,n)}^{\text{con}})=\chi_{\1}(I_{(k,n)}^{\text{con}})$, obeys the following recursion (for $n\geq 2k+1$) 
\[
d(\Sc_{(k,n)}^{\text{con}})=d(\Sc_{(k,n-(k+1))}^{\text{con}}),
\] 

with base cases 
\begin{align*}
    d(\Sc_{(k,k)}^{\text{con}})&=1\\
    d(\Sc_{(k,k+1)}^{\text{con}})&=-1, \text{ and}\\
    d(\Sc_{(k,n)}^{\text{con}})&=0, \text{ for all }k+1<n\leq 2k,
\end{align*}
from which the results follows.

Let us first check the base cases:
\begin{itemize}
    \item{\fbox{$n=k$}: The ideal associated to this system is $I_{(k,k)}^{\text{con}}=\langle x_1\cdots x_k\rangle$. Clearly $\beta_{0,\1}=1$ and $\beta_{i,\1}=0$, for all $i>0$, hence, $\chi_{\1}(I_{(k,k)}^{\text{con}})=1$.
    }
    
    \item{\fbox{$n=k+1$}: The ideal of a binary consecutive linear $k$-out-of-$k+1$ system is $I_{(k,k+1)}^{\text{con}}=\langle x_1\cdots x_{k},x_2\cdots x_{k+1}\rangle$, which only has two generators. Clearly, $\beta_{1,\1}=1$ and $\beta_{i,\1}=0, \forall i\neq 1$, hence $\chi_{\1}(I_{(k,k+1)}^{\text{con}})=(-1)^1=-1$.
    }

    \item{\fbox{$k+1<n\leq 2k$}: Consider the MVT of $I_{(k,n)}^{\text{con}}$, with the generators sorted lexicographically in each node, and take always the last generator as a pivot. The node at position $3$ has the generators of the ideal $I_{(k,n-1)}^{\text{con}}$, in which the variable $x_n$ does not appear. Hence, this subtree does not contribute to $\chi_{\1}(I_{(k,n)}^{\text{con}})$.

    The node at position $2$ is generated by $x_{n-k}\cdot x_{n-k+1}\cdots x_n$, in which the first $n-k-1$ variables are not present. Hence, this subtree does not contribute to $\chi_{\1}(I_{(k,n)}^{\text{con}})$ either. To see that MVT$_{2,1}=\langle x_{n-k}\cdot x_{n-k+1}\cdots x_n\rangle$ observe that $\lcm(x_{n-k}\cdots x_{n-1},x_{n-k+1}\cdots x_n))=x_{n-k}\cdot x_{n-k+1}\cdots x_n$
    and that the least common multiple of any other generator of $I_{(k,n)}^{\text{con}}$ with $x_{n-k+1}\cdots x_n$ also contains $x_{n-k}$ and, hence, is divided by $x_{n-k}\cdot x_{n-k+1}\cdots x_n$.
    
    Since neither of $MVT_{3,0}(I_{(k,n)}^{\text{con}})$ and $MVT_{2,1}(I_{(k,n)}^{\text{con}})$ contribute to $\chi_{\1}(I_{(k,n)}^{\text{con}})$, we conclude that $\chi_{\1}(I_{(k,n)}^{\text{con}})=0$.
 }
\end{itemize}

\begin{figure}[h]
\centering
\resizebox{\textwidth}{!}{%
\begin{tikzpicture}[
  every node/.style={anchor=west, inner sep=2pt},
  flecha/.style={-{Stealth[length=2mm]}, gray!70!black, shorten <=1pt, shorten >=2pt}
]
\node (a) at (0,0)
  {$(1,0)\quad x_1x_2\cdots x_k,\ x_2x_3\cdots x_{k+1},\ \dots,\ x_{n-k+1}x_{n-k+2}\cdots x_n$};
\node (b) at (3.5,-1.4)
  {$(3,0)\quad x_1x_2\cdots x_k,\ x_2x_3\cdots x_{k+1},\ \dots,\ x_{n-k}x_{n-k+1}\cdots x_{n-1}$};
\node (c) at (0.8,-2.8)
  {$(2,1)\quad (x_{n-k+1}\cdots x_n)\,(x_1\cdots x_k,\ x_2\cdots x_{k+1},\ \dots,\ x_{n-2k}\cdots x_{n-k-1},\ x_{n-k})$};
\node (d) at (4.3,-4.2)
  {$(5,1)\quad (x_{n-k+1}\cdots x_n)\,(x_1\cdots x_k,\ x_2\cdots x_{k+1},\ \dots,\ x_{n-2k}\cdots x_{n-k-1})$};
\node (e) at (1.6,-5.6)
  {$(4,2)\quad (x_{n-k}x_{n-k+1}\cdots x_n)\,(x_1\cdots x_k,\ x_2\cdots x_{k+1},\ \dots,\ x_{n-2k}\cdots x_{n-k-1})$};

\draw[flecha] ([xshift=6mm]a.south west) -- (b.west);
\draw[flecha] ([xshift=6mm]a.south west) -- ([xshift=4mm]c.north west);
\draw[flecha] ([xshift=6mm]c.south west) -- (d.west);
\draw[flecha] ([xshift=6mm]c.south west) -- ([xshift=4mm]e.north west);
\end{tikzpicture}%
}
\caption{Partial MVT of a binary consecutive linear $k$-out-of-$n$ system.}
\label{fig:MVT_consKN}
\end{figure}
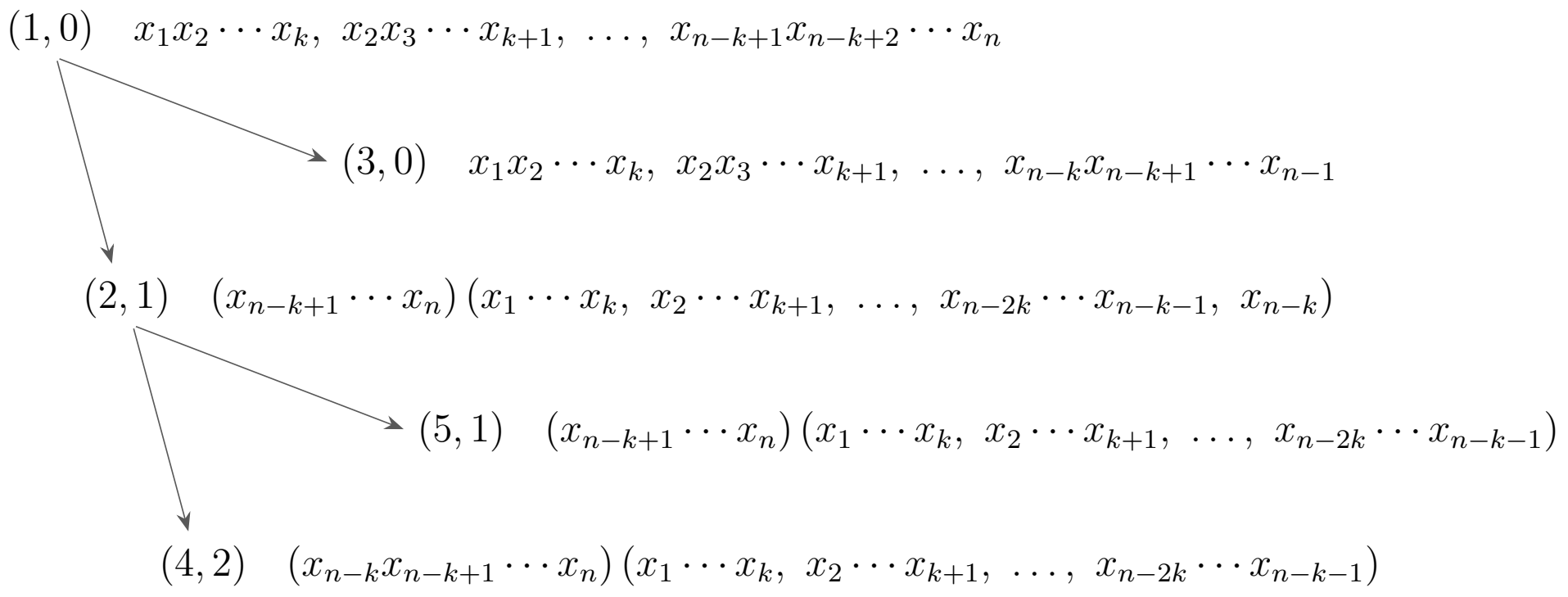

%

Finally, to verify that the recursion $d(\Sc_{(k,n)}^{\text{con}})=d(\Sc_{(k,n-(k+1))}^{\text{con}})$ is satisfied, let us consider $MVT(I_{(k,n)}^{\text{con}})$ with generators sorted lexicographically, taking the last generator as a pivot. See Figure \ref{fig:MVT_consKN}\footnote{Observe that in this figure, the nodes in positions $2,4$ and $5$ are written, for the sake of clarity, in the form $\mu\cdot(\nu_1,\dots,\nu_l)$ for some $l$, which means $\mu\nu_1,\dots,\mu\nu_l$.}. The subtree from node $MVT_{3,0}$ corresponds to the MVT of $I_{(k,n-1)}^{\text{con}}$ which does not contribute to $\chi_{\1}(I_{(k,n)}^{\text{con}})$.

On the other hand, the generators of node $MVT_{2,1}$ are all of them of the form $x_{n-k+1}\cdots x_n\cdot\mu$, where  $$\mu\in\left\lbrace x_1\cdots x_k, x_2\cdots x_{k+1}, \dots, x_{n-2k}\cdots x_{n-(k+1)}\right\rbrace$$ or $\mu=x_{n-k}$. Taking now $x_{n-k+1}\cdots x_n\cdot x_{n-k}$ as pivot in this node, we obtain the following:
\begin{itemize}
    \item{In node $MVT_{5,1}$ the variable $x_{n-k}$ is not present. Then the subtree hanging from $MVT_{5,1}$ does not contribute to $\chi_{\1}(I_{(k,n)}^{\text{con}})$.
    }
    \item{The node $MVT_{4,2}$ is generated by all elements of the form $x_{n-k}\cdots x_{n}\cdot\mu$, where $$\mu\in\left\lbrace x_1\cdots x_k, x_2\cdots x_{k+1},\dots, x_{n-2k}\cdots x_{n-(k+1)}\right\rbrace$$. Taking the common factor $x_{n-k}\cdots x_n$ in these generators, we obtain the generators of $I_{(k,n-(k+1))}^{\text{con}}$. Hence the contribution of the subtree hanging from node $MVT_{4,2}$ to $\chi_{\1}(I_{(k,n)}^{\text{con}})$ is exactly $\chi_{x_1\cdots x_{n-(k+1)}}(I_{(k,n-(k+1))}^{\text{con}})$, and the recursion is proven. 
    }
\end{itemize}

\end{proof}


\begin{remark}
Observe that for $k=2$, Proposition \ref{prop:BettiConsecBinKN} gives the Euler characteristic of the independence complex of a path up to the sign $(-1)^n$ (see \cite{K99}). In general, the consecutive binary system ideal is the $k$-path ideal of a path graph \cite{BHO11}.
\end{remark}
\subsection{Multi-state $k$-out-of-$n$ and consecutive linear $k$-out-of-$n$ systems.}\label{sec:multi_k-out-of-n}
The definition of multi-state $k$-out-of-$n$ and consecutive linear multi-state $k$-out-of-$n$ systems is not as direct as in the binary case, and several different definitions have been given in the literature. In this section we study the simple and generalized versions of these kind of systems.

\subsubsection{Simple multi-state $k$-out-of-$n$ systems.}
The first definition of multi-state $k$-out-of-$n$ systems is the following, due to El-Newehi, Proshan and Sethuraman \cite{NPS78}:

\begin{definition}[El-Neweihi et al., $1978$] \label{Def:KN_EN_78}
	A system is a \emph{multi-state $k$-out-of-$n$ system} if its structure function satisfies
	\begin{equation}\label{eq:KN_EN_78}
	 \phi(\xb)=x_{(n-k+1)}
	 \end{equation}
	 where $x_{(1)}\leq x_{(2)}\leq\dots\leq x_{(n)}$ is a non decreasing arrangement of $x_1,\dots, x_n$. 
\end{definition}

This was generalized by Boedigheimer and Kapur \cite{BK94}, using the notion of lower and upper boundary points.

\begin{definition}
	A state vector $\xb$ is a \textit{lower boundary point {\rm (l.b.p.)} to level $j$} if and only if $\phi(\xb)\geq j$ and $\yb<\xb$ implies $\phi(\yb)<j$, for $j=1,\dots,M$.
	An \textit{upper boundary point {\rm (u.b.p)} to level $j$} is a state vector $\xb$ such that $\phi(\xb)\leq j$ and $\yb>\xb$ implies $\phi(\yb)>j$, for $j=0,\dots,M-1$.
\end{definition}

\begin{definition}[Boedigheimer and Kapur, $1994$]\label{Def:KN_BK_94}
	A function $\phi$ is the structure function of a simple multi-state $k$-out-of-$n:G$ iff, $\phi$ has ${n_j \choose k}$ lower boundary points to level $j$ $(j=1,\dots, M)$ and ${n_j \choose k-1}$ upper boundary points to level $j$ $(j=0,\dots, M-1)$, where $n_j\leq n$ be the number of components such that $M_i\geq j$ for each $j\in\{0,\dots,M\}$.
\end{definition}

In what follows, we use Definition \ref{Def:KN_BK_94} and denote this kind of systems as $\Sc_{(k,n),M}$. The $j$-reliability ideals of simple multi-state $k$-out-of-$n$ systems can be described as follows, see \cite{PSW21}.

\begin{proposition}\label{prop:K_N_alg}
Let $\Sc_{(k,n),M}$ be a simple multi-state $k$-out-of-$n$ system with levels $\{0,\dots,M\}$ and such that each component $i$ has a set of increasing levels of performance given by $\{0,\dots, M_i\}$. Let $n_j\leq n$ be the number of components such that $M_i\geq j$ for each $j\in\{0,\dots,M\}$ (without loss of generality, we assume that these are components $1,\dots,n_j$). For every $j\in\{1,\dots,M\}$ the $j$-reliability ideal of $\Sc_{(k,n),M}$, is 
 \[
 I_{(k,n_j),j}=\left\langle \prod_{\substack{\sigma\subseteq\lbrace 1,\dots, n_j\rbrace \\ \vert\sigma\vert=k}} x_i^j\ \vert\, i\in\sigma \right\rangle.
 \] 
\end{proposition}



The domination of these systems is usually equal to $0$ except for some special cases, in which it is equivalent to that of some binary $k$-out-of-$n$ systems.

\begin{proposition}\label{prop:simpleMS-KN}
    Let $\phi$ be the structure function of a simple multi-state $k$-out-of-$n:G$ system $\Sc_{(k,n),M}$. For every level $j$ with $n_j\geq k$,
    \begin{align*}
    \overline d_j(\Sc_{(k,n),M})&=(-1)^{n_j-k}{{n_j-1}\choose{k-1}}.
    \end{align*}
\end{proposition}
\begin{proof}
The ideal $I_{(k,n_j),j}$ is obtained from the binary ideal $I_{(k,n_j)}$ by the substitution $x_i\rightarrow x_i^j$, hence the maximal element in $L_{I_{(k,n_j),j}}$ is $x_1^j\cdots x_{n_j}^j$. Then, by Proposition \ref{prop:dom_binaryKN}, $\overline d_j(\Sc_{(k,n),j})=(-1)^{n_j-k}{{n_j-1}\choose{k-1}}$.
\end{proof}

\begin{remark}
    Observe that for these systems, those variables $x_i$ such that $i>n_j$ are not relevant for level $j$, hence the system is not coherent at level $j$ unless $n_j=n$. To avoid this, we can consider that at level $j$, the system consists only on the components $\{c_1,\dots,c_{n_j}\}$ and consider $\overline d_j(\Sc_{(k,n),j})$.
\end{remark}

\begin{example}[Example 3.10 in \cite{PSW21}]\label{ex:simpleMultiKN}
Consider a multi-state system $\Sc$ such that $S_1 = \{0, 1, 2, 3, 4\}, S_2 = \{0, 1, 2, 3\}, S_3 = S_4 =
\{0, 1, 2\}$ and $S_5 = \{0, 1\}$ and let $\phi(x) = x_{(4)}$ (see Definition \ref{Def:KN_EN_78}). The upper and lower boundary points of this system are listed in Table \ref{table:exampleSimpleMultiKN}, and from them we see that, according to Definition \ref{Def:KN_BK_94}, the system behaves as a $2$-out-of-$n_j$ system for levels $j=1,2,3$, where $n_1 = 5, n_2 = 4, n_3 = 2$.

\begin{table}[htpb]
  \centering
  \caption{Upper and lower boundary points for the system in Example \ref{ex:simpleMultiKN}}
  \label{table:exampleSimpleMultiKN}
  \begin{tabularx}{\linewidth}{l >{\hsize=1.3\hsize\raggedright\arraybackslash}X >{\hsize=0.7\hsize\raggedright\arraybackslash}X}
    \toprule
    Level & Lower boundary points & Upper boundary points \\
    \midrule
    $0$ & & $(0, 0, 0, 0, 1), (0, 0, 0, 2, 0), $ \newline $(0, 0, 2, 0, 0),(0, 3, 0, 0, 0), $ \newline $(4, 0, 0, 0, 0)$\\
    \addlinespace
    $1$ & $(0, 0, 0, 1, 1), (0, 0, 1, 0, 1), (0, 1, 0, 0, 1), $ \newline $(1, 0, 0, 0, 1),(0, 0, 1, 1, 0),(0, 1, 0, 1, 0), $ \newline $(1, 0, 0, 1, 0), (0, 1, 1, 0, 0)(1, 0, 1, 0, 0),$ \newline $ (1, 1, 0, 0, 0)$ 
      & $(1, 1, 1, 2, 1), (1, 1, 2, 1, 1),$ \newline $(1, 3, 1, 1, 1),(4, 1, 1, 1, 1)$ \\
    \addlinespace
    $2$ & $(0, 0, 2, 2, 0), (0, 2, 0, 2, 0), (2, 0, 0, 2, 0), $ \newline $(0, 2, 2, 0, 0),(2, 0, 2, 0, 0), (2, 2, 0, 0, 0)$
      & $(2, 3, 2, 2, 1), (4, 2, 2, 2, 1)$ \\
    \addlinespace
    $3$ & $(3, 3, 0, 0, 0)$ & \\
    \bottomrule
  \end{tabularx}
\end{table}

The $j$-reliability ideals of this system are then generated by the monomials corresponding to the lower boundary points to each of the levels, specifically 
\begin{align*}
J_1(\Sc)=&\langle x_1x_2,x_1x_3,x_1x_4,x_1x_5,x_2x_3,x_2x_4,x_2x_5,x_3x_4,x_3x_5,x_4x_5\rangle,\\
J_2(\Sc)=&\langle x_1^2x_2^2,x_1^2x_3^2,x_1^2x_4^2,x_2^2x_3^2,x_2^2x_4^2,x_3^2x_4^2 \rangle,\\
J_3(\Sc)=&\langle x_1^3x_2^3\rangle,
\end{align*}
and it is easy to see that $d_j(\Sc)=\chi_{(4,3,2,2,1)}(I_{(2,n_j),j})=0$ for $j=1,2,3$. Observe that the system is coherent for $j=1$ but it is not coherent for $j=2,3$. The level dominations for this system are $\overline d_1(\Sc)=-4, \overline d_2(\Sc)=3,\overline d_3(\Sc)=1$.
\end{example}


\subsubsection{Generalized multi-state $k$-out-of-$n$ systems.}
A more general definition of multi-state $k$-out-of-$n$ systems was given in \cite{ZLW00,ZT06}. Both :G (for \emph{good}) and :F (for \emph{failure}) types of systems are defined in those papers; we restrict here to :G systems. In these subsection we consider that the maximum level of performance $M$ of the system, equals the maximum performance level $M_i$ of all its components.
\begin{definition}
An $n$-component system is called a generalized multi-state $k$-out-of-$n$:G system, denoted by $\Sc_{(K,n),M}$, if $\phi(x)\geq j, 1\leq j\leq M$ whenever there exists an integer value $l$ $(j\leq l\leq M)$ such that at least $k_l$ components are in state $l$ or above.
\end{definition}
We can describe these systems in the following way. Consider a sequence of integers $K=(k_1,\dots,k_M)$, the system is in level $j$ or above if at least $k_j$ components are in level $j$ or above, at least $k_{j+1}$ components are in level $j+1$ or above ... or at least $k_M$ components are in level $M$. Observe that these are always coherent systems.

The $j$-reliability ideal of a generalized multi-state $k$-out-of-$n$ system $\Sc_{(K,n),M}$ defined by the sequence $K=(k_1,\dots,k_M)$ is the following (see \cite{PSW21}).
\[
I_{(K,n),j}=\sum_{i=j}^M I_{(k_i,n),i}
\]
where the ideals $I_{(k_i,n),i}$ are as in Proposition \ref{prop:K_N_alg}.

\begin{proposition}\label{prop:multistateKN}
Let $S_{(K,n),M}$ be as in Definition~7.5, with $M_i=M$ for all $i$ and
$1\le k_l\le n$. For every level $j$,
\[
  \overline d_j(S_{(K,n),M})=
  \begin{cases}
   (-1)^{n-k_j} {{n-1}\choose {k_j-1}} & \text{if } k_j\le k_l \text{ for all } l>j,\\[8pt]
   0 & \text{otherwise.}
  \end{cases}
\]
\end{proposition}
 
\begin{proof}
Every generator of $I_{(k_l,n),l}$ with $l>j$ ($k_l$ variables to the power $l$)
is divisible by some generator of $I_{(k_{l'},n),l'}$, $j\le l'<l$, iff
$k_{l'}\le k_l$.
 
If $k_j\le k_l$ for all $l>j$, all generators of the higher levels are
redundant. Then $I_{(K,n),j}=I_{(k_j,n),j}$ and $\nu_j=x_1^j\cdots x_n^j$, and the
substitution $x_i\mapsto x_i^j$ together with Proposition \ref{prop:dom_binaryKN} gives the
value.
 
Otherwise let $l_0$ be the largest $l>j$ with $k_l<k_{l'}$ for all
$j\le l'<l$. Its generators are minimal, and no higher level contributes
minimal generators, so $\nu_j=x_1^{l_0}\cdots x_n^{l_0}$ with $l_0>j$. Then
$\nu_j/x_1\cdots x_n$ is in the ideal $I_{(K,n),j}$, hence $\Delta K(I_j,\nu_j)$ is the full simplex, whose reduced Euler characteristic is $0$, and then $\overline d_j=0$ by Theorem \ref{th:dominationKoszul}.
\end{proof}

\begin{remark}
Observe that for these systems $\overline d_M$ is always nonzero; all level dominations are nonzero iff $K$
is non-decreasing; and only $\overline d_M$ is nonzero iff $k_M<k_l$ for all $l<M$.
\end{remark}

\begin{example}
 For $n=4$, $M=4$, $K=(4,1,2,3)$ (all values
confirmed by brute force):
\begin{center}
\renewcommand{\arraystretch}{1.3} 
\begin{tabular}{clll}
\toprule
$j$ & minimal generators of $I_{(K,4),j}$ & $\nu_j$ & $\overline d_j$ \\
\midrule
$4$ & $I_{(3,4),4}$ & $(4,4,4,4)$ & $(-1)^{1}\binom32=-3$ \\
\addlinespace
$3$ & $I_{(2,4),3}$ \mbox{ (level 4 is redundant, $k_3\le k_4$)} & $(3,3,3,3)$ & $(-1)^{2}\binom31=3$ \\
\addlinespace
$2$ & $\langle x_1^2,\dots,x_4^2\rangle$\mbox{ (levels 3, 4 are redundant}) & $(2,2,2,2)$ & $(-1)^{3}\binom30=-1$ \\
\addlinespace
$1$ & $\langle x_1x_2x_3x_4\rangle+\langle x_1^2,\dots,x_4^2\rangle$ & $(2,2,2,2)$ & $0$ \\
\bottomrule
\end{tabular}
\end{center}
\end{example}

\subsubsection{Simple multi-state consecutive $k$-out-of-$n$ systems.}

A simple multi-state consecutive $k$-out-of-$n$ system, denoted by $\Sc_{(k,n),M}^{\text{con}}$, is a multi-state system with $n$ components that works at level $j$ if at least $k$ consecutive components are working at level $j$.

The $j$-reliability ideal of such a system is of the form \[I_{(k,n),j}^{\text{con}}=\langle x_1^j\dots x_k^j,x_2^j\dots x_{k+1}^j,\dots,x_{n-k+1}^j\dots x_n^j \rangle.\]

\begin{proposition}\label{prop:BettiConsecSMSKN}
    Let $\Sc_{(k,n),M}^{\text{con}}$ be a simple multi-state consecutive $k$-out-of-$n$:$G$ system, with $0<k\leq n$, and $I_{(k,n),j}^{\text{con}}$ its $j$-reliability ideal. Then, if $j$ is the maximum level of performance of all the components of the system (i.e. $M_i=j$ for all $i$), the domination of this system is given by
	 
		\centering
		$d_j(\Sc_{(k,n),M}^{\text{con}}) = \begin{cases}
        1 & \text{if } n\equiv k\mod(k+1),\\
        -1  & \text{if } n\equiv 0\mod(k+1),\\
        0  & \text{otherwise}.
        \end{cases}$
\end{proposition}

\begin{proof}
    The proof is essentially the same as in Proposition \ref{prop:BettiConsecBinKN}.
\end{proof}

\subsubsection{Generalized multi-state consecutive-$k$-out-of-$n$ systems.}\label{sec:generalizedConsecutiveKN}
Generalized multi-state consecutive $k$-out-of-$n$ :F and :G were defined in \cite{HZF03} as follows (we consider here only the definition of :G systems).
\begin{definition}
A system with a structure function $\phi$ such that $\phi(\xb) \geq j$, $( j=1, 2,\dots, M)$ if at least $k_l$ consecutive components are in state $l$ or above for all $l$, $(1\leq l \leq j)$ is called a multi-state consecutive-$k$-out-of-$n$:G system.
\end{definition}

We will denote them as $\Sc_{(K,n),M}^{\text{con}}$.

These systems have a $k_l$-out-of-$n$ structure at each of the levels, which represents the condition of having at least $k_l$ consecutive components performing at level $l$ or higher. The whole system must satisfy all the conditions at the same time. The ideals corresponding to each of the level conditions are simple multi-state consecutive $k$-out-of-$n$ systems, denoted by $I_{(k_l,n),l}^{\text{con}}$, hence the ideal of a generalized multi-state consecutive $k$-out-of-$n$  system is given by
\[
I_{(K,n),j}^{\text{con}}=\bigcap_{1\leq l\leq j}I_{(k_l,n),l}^{\text{con}},
\]
for $K=(k_1,\dots,k_j)$.

The domination of such a system is given by the domination of the system corresponding to the last level, which equals the domination of a binary consecutive $k_j$-out-of-$n$ system.
\begin{proposition}

\[
\overline d_j(\Sc_{(K,n),M}^{\text{con}})=\overline d_j(\Sc_{(k_j,n),M}^{\text{con}})=\overline d_j(\Sc_{(k_j,n)}^{\text{con}})= \begin{cases}
    1, & \text{if }  n\equiv k_j \bmod (k_j+1),\\
    -1 , & \text{if } n\equiv 0 \bmod (k_j+1),\\
    0, & \mbox{otherwise.}
\end{cases}
\]
\end{proposition}

\begin{proof}
We prove this result by showing that the Koszul complex of $d(\Sc_{(K,n),M}^{\text{con}})$ and $d(\Sc_{(k_j,n),M})$ at $\mu=(j,\dots,j)$ coincide, and therefore both dominations are the same. The last equality is Proposition \ref{prop:BettiConsecSMSKN}.

Recall that $G(I_{(k_j,n),j}^{\text{con}})$ consists on all possible products of $k_j$ consecutive variables, each of them raised to the power $j$. On the other hand, each of minimal monomial generators of $I_{(K,n)}^{\text{con}}$ consists on a product of $k_j$ consecutive variables, each raised to the power $j$, times some other variables, raised to a strictly smaller exponent. Note that in this last ideal, each possible product of $k_j$ consecutive variable appears in at least one minimal monomial generator of $I_{(K,n)}^{\text{con}}$. 

Consider now
\[
T=\{ \sigma\subseteq\{1,\dots,n\}\mid \vert\sigma\vert=n-k_j \mbox{ and } \overline\sigma \mbox{ consists of } k_j\mbox{ consecutive indices}\}.
\]

We claim that $T$ is the set of facets of both $\Delta K(I_{(k_j,n),j})$ and $\Delta K(I_{(K,n)})$, and therefore both complexes are equal.

Let $\sigma\in T$, the monomial $\xb^{\mu/\sigma}$ is by definition, the product of the $n-k_j$ variables in $\sigma$ raised to the power $j-1$ times the variables in $\overline{\sigma}$ each raised to the power of $j$, these are exactly $k_j$ consecutive variables, hence there is one and only one generator in $I_{(k_j,n)}$ that divides $\xb^{\mu/\sigma}$, hence $\sigma\in \Delta K(I_{(k_j,n),j})$. Also, there is at least one generator of $I_{(K,n)}$ that divides $\mu/\sigma$, hence $\sigma\in \Delta K(I_{(K,n)})$.

On the other hand, let $\tau\subseteq\{1,\dots, n\}$ such that there is no $\sigma\in T$ that satisfies $\tau\subseteq\sigma$. Observe that for such a set $\tau$, we have that $\overline{\tau}$ does not contain any set of $k_j$ consecutive indices, hence $\mu/\tau$ is not divisible by any set of $k_j$ consecutive variables raised to the power $j$, and therefore there is no minimal generator of $I_{(k_j,n),j}$ or $I_{(K,n)}$ that divides $\xb^{\mu/\tau}$. Hence, by definition, $\tau\notin \Delta K(I_{(k_j,n),j})$ and $\tau \notin \Delta K(I_{(K,n)})$.

\end{proof}

\section{Conclusions}\label{sec:conclusions}
In this paper we described algebraic and topological tools to analyze and compute the signed domination function of binary and multi-state systems. The correspondence between such systems and monomial ideals opened the possibility of these tools to study the reliability of such systems, and was used by the authors and others in several papers during the last decade. In this paper, we went a step further, establishing in our main result a correspondence between the signed domination function and the Euler characteristic of combinatorial objects associated to the ideal corresponding to the system, namely chain complexes (resolutions in this case) and simplicial complexes. Based on this, we propose a new method to compute the signed domination function and domination of binary and multi-state systems. We apply this method to obtain the signed domination function and signed domination of several systems, that include component-based dominant systems and variants of $k$-out-of-$n$ systems, both in the binary and in the multi-state case.

\section*{Acknowledgements}
This work has been partially supported by and grant PID2024-157733NBI00 funded by MCIN/AEI/10.13039/501100011033/FEDER EU.

\section{Declaration of generative AI and AI-assisted technologies in the manuscript preparation process.}
During the preparation of this work the authors used the AI tool Claude (model Opus 5.5) to check the paper in terms of notation, correctness of the statements and proofs and grammar and spelling checks. No new results or other significant additions were included as a result of the interaction with the AI tool. After using this tool/service, the authors reviewed and edited the content as needed and take full responsibility for the content of the published article.

\bibliographystyle{plain}
\bibliography{dom_bibliography.bib}

@article{AB84,
  author  = {Agrawal, A. and Barlow, R. E.},
  title   = {A survey of network reliability and domination theory},
  journal = {Operations Research},
  volume  = {32},
  number  = {3},
  pages   = {478--492},
  year    = {1984}
}

@article{A17,
  title={Minimal resolutions of dominant and semidominant ideals},
  author={Alesandroni, Guillermo},
  journal={Journal of Pure and Applied Algebra},
  volume={221},
  number={4},
  pages={780--798},
  year={2017},
  publisher={Elsevier}
}

@article{AFL22,
  title={Polarizations of powers of graded maximal ideals},
  author={Almousa, A. and Fl{\o}ystad, G. and Lohne, H.},
  journal={Journal of Pure and Applied Algebra},
  volume={226},
  number={5},
  pages={106924},
  year={2022},
  publisher={Elsevier}
}

@book{AM69,
  title={Introduction to Commutative Algebra},
  author={Atiyah, M.F. and MacDonald, I.G.},
  year={1969},
  publisher={Addison-Wesley}
}

@article{AV22,
  title={Polarizations and hook partitions},
  author={Almousa, A. and VandeBogert, K.},
  journal={Journal of Pure and Applied Algebra},
  volume={226},
  number={9},
  pages={107056},
  year={2022},
  publisher={Elsevier}
}

@article{B86,
  title={Computational complexity of network reliability analysis: An overview},
  author={Ball, M. O.},
  journal={IEEE Transactions on Reliability},
  volume={35},
  number={3},
  pages={230--239},
  year={1986},
  publisher={IEEE}
}

@techreport{B82,
  title={Set Theoretic Signed Domination for Coherent Systems.},
  author={Barlow, R. E},
  year={1982},
  institution={University of California},
  address={Berkeley}
}

@article{BS98,
title = {Cellular resolutions of monomial modules},
author={Bayer, D. and Sturmfels, B.},
pages = {123--140},
volume = {1998},
number = {502},
journal = {Journal für die reine und angewandte Mathematik},
doi = {doi:10.1515/crll.1998.083},
year = {1998}
}

@inproceedings{BHS19,
  title={Monomial resolutions for efficient computation of simplicial homology},
  author={Bigatti, A. M. and Heras, J. and {S{\'a}enz-de-Cabez{\'o}n}, E.},
  booktitle={Proceedings of the 2019 International Symposium on Symbolic and Algebraic Computation},
  pages={50--57},
  year={2019}
}

@article{BPS21,
  title={A {C++} class for multi-state algebraic reliability computations},
  author={Bigatti, A. M. and Pascual-Ortigosa, P. and S{\'a}enz-de-Cabez{\'o}n, E.},
  journal={Reliability Engineering \& System Safety},
  volume={213},
  pages={107751},
  year={2021},
  publisher={Elsevier}
}

@article{BS82,
  title={A decomposition for multistate monotone systems},
  author={Block, H. W. and Savits, T. H.},
  journal={Journal of Applied Probability},
  volume={19},
  number={2},
  pages={391--402},
  year={1982},
  publisher={Cambridge University Press}
}

@article{BK94,
  title={Customer-driven reliability models for multistate coherent systems},
  author={Boedigheimer, R. A and Kapur, K. C.},
  journal={IEEE Transactions on Reliability},
  volume={43},
  number={1},
  pages={46--50},
  year={1994},
  publisher={IEEE}
}

@article{B01,
  author  = {Boland, P. J.},
  title   = {Signatures of indirect majority systems},
  journal = {Journal of Applied Probability},
  volume  = {38},
  number  = {2},
  pages   = {597--603},
  year    = {2001}
}

@article{BHO11,
  title={Path ideals of rooted trees and their graded Betti numbers},
  author={Bouchat, R. R. and H\'a, H. T. and O'Keefe, A.},
  journal={Journal of Combinatorial Theory, Series A},
  volume={118},
  number={8},
  pages={2411--2425},
  year={2011},
  publisher={Elsevier}
}

@article{CE95,
  title={Resolutions obtained by iterated mapping cones},
  author={Charalambous, H. and Evans, E. G.},
  journal={Journal of Algebra},
  volume={176},
  number={3},
  pages={750--754},
  year={1995},
  publisher={Academic Press}
}

@article{C10,
  title={Poset resolutions and lattice-linear monomial ideals},
  author={Clark, T. BP},
  journal={Journal of Algebra},
  volume={323},
  number={4},
  pages={899--919},
  year={2010},
  publisher={Elsevier}
}

@book{CLO05,
  title={Using algebraic geometry},
  author={Cox, D. and Little, J. and O'shea, D.},
  volume={185},
  year={2005},
  publisher={Springer}
}

@book{D03,
  title={Improved Bonferroni inequalities via abstract tubes: inequalities and identities of inclusion-exclusion type},
  author={Dohmen, K.},
  year={2003},
  publisher={Springer}
}

@article{EK90,
  title={Minimal resolutions of some monomial ideals},
  author={Eliahou, S. and Kervaire, M.},
  journal={Journal of Algebra},
  volume={129},
  number={1},
  pages={1--25},
  year={1990},
  publisher={Academic Press}
}

@article{NPS78,
  title={Multistate coherent systems},
  author={El-Neweihi, E. and Proschan, F. and Sethuraman, J.},
  journal={Journal of Applied Probability},
  volume={15},
  number={4},
  pages={675--688},
  year={1978},
  publisher={Cambridge University Press}
}

@incollection{F05,
  title={Monomial ideals via square-free monomial ideals},
  author={Faridi, S.},
  booktitle={Commutative Algebra: Geometric, Homological, Combinatorial and Computational Aspects},
  editor= {Corso, A. and Gimenez, P. and Vaz Pinto, M. and Zarzuela, S.},
  pages={107--136},
  year={2005},
  publisher={CRC Press}
}

@article{FMS16,
  title={LCM lattices supporting pure resolutions},
  author={Francisco, C. and Mermin, J. and Schweig, J.},
  journal={Proceedings of the American Mathematical Society},
  volume={144},
  number={6},
  pages={2315--2325},
  year={2016}
}

@ARTICLE{GPW99,
	author = {Gasharov, V. and Peeva, I. and Welker, V.},
	title = {The LCM-lattice in monomial resolutions},
	year = {1999},
	journal = {Mathematical Research Letters},
	volume = {6},
	number = {5-6},
	pages = {521--532},
	doi = {10.4310/MRL.1999.v6.n5.a5},
}

@article{GW04,
  title={Monomial ideals and the {S}carf complex for coherent systems in reliability theory},
  author={Giglio, B. and Wynn, H. P},
  year={2004},
  journal = {Annals of Statistics},
  volume = {32},
  number = {3},
  pages = {1289--1311}
}

@book{H04,
  title={Algebraic Topology},
  author={Hatcher, A.},
  year={2001},
  publisher={Cambridge University Press}
}

@book{HH11,
  title={Monomial ideals},
  author={Herzog, J. and Hibi, T.},
  year={2011},
  publisher={Springer}
}

@article{H77,
  title={Cohen-{M}acaulay rings, combinatorics, and simplicial complexes, in" Ring Theory II"},
  author={Hochster, M.},
  journal={Lect. Notes in Pure Appl. Math.},
  number={26},
  pages={171--223},
  year={1977},
  publisher={Dekker}
}

@article{HZF03,
  title={Multi-state consecutive-k-out-of-n systems},
  author={Huang, J. and Zuo, M. J and Fang, Z.},
  journal={IIE Transactions},
  volume={35},
  number={6},
  pages={527--534},
  year={2003},
  publisher={Taylor \& Francis}
}

@article{HZ04,
  title={Dominant multi-state systems},
  author={Huang, J. and Zuo, M. J.},
  journal={IEEE transactions on reliability},
  volume={53},
  number={3},
  pages={362--368},
  year={2004},
  publisher={IEEE}
}

@article{H89,
  author  = {Huseby, A. B.},
  title   = {Domination theory and the {C}rapo $\beta$-invariant},
  journal = {Networks},
  volume  = {19},
  number  = {1},
  pages   = {135--149},
  year    = {1989}
}

@article{H25,
      title={Domination and Multistate Systems},
  author={Huseby, A. B.},
  journal={arXiv preprint arXiv:2502.17737},
  year={2025}
}

@article{K99,
  title={Complexes of directed trees},
  author={Kozlov, D. N},
  journal={Journal of Combinatorial Theory, Series A},
  volume={88},
  number={1},
  pages={112--122},
  year={1999},
  publisher={Elsevier}
}

@book{KZ03,
  title={Optimal reliability modeling: principles and applications},
  author={Kuo, W. and Zuo, M. J},
  year={2003},
  publisher={John Wiley \& Sons}
}

@article{M13,
  title={Finite atomic lattices and resolutions of monomial ideals},
  author={Mapes, S.},
  journal={Journal of Algebra},
  volume={379},
  pages={259--276},
  year={2013},
  publisher={Elsevier}
}

@article{MM11,
  author  = {Marichal, J.-L. and Mathonet, P.},
  title   = {Extensions of system signatures to dependent lifetimes: Explicit expressions and interpretations},
  journal = {Journal of Multivariate Analysis},
  volume  = {102},
  number  = {5},
  pages   = {931--936},
  year    = {2011}
}

@article{MMVV19,
  title={Depth and regularity of monomial ideals via polarization and combinatorial optimization},
  author={Mart{\'\i}nez-Bernal, J. and Morey, S. and Villarreal, R. H and Vivares, C. E},
  journal={Acta Mathematica Vietnamica},
  volume={44},
  number={1},
  pages={243--268},
  year={2019},
  publisher={Springer}
}

@article{M08,
  title={Borel-plus-powers monomial ideals},
  author={Murai, S.},
  journal={Journal of Pure and Applied Algebra},
  volume={212},
  number={6},
  pages={1321--1336},
  year={2008},
  publisher={Elsevier}
}

@article{IMPSW24,
  title={Stable coherent systems},
  author={Iglesias, R. and Mohammadi, F. and Pascual-Ortigosa, P. and S{\'a}enz-de-Cabez{\'o}n, E. and Wynn, H. P},
  journal={Applicable Algebra in Engineering, Communication and Computing},
  pages={1--27},
  year={2024},
  publisher={Springer}
}

@book{MS05,
  title={Combinatorial commutative algebra},
  author={Miller, E. and Sturmfels, B.},
  volume={227},
  year={2005},
  publisher={Springer Verlag}
}

@article{MPSW20,
  title={Polarization and depolarization of monomial ideals with application to multi-state system reliability},
  author={Mohammadi, F. and Pascual-Ortigosa, P. and S{\'a}enz-de-Cabez{\'o}n, E. and Wynn, H. P},
  journal={Journal of Algebraic Combinatorics},
  volume={51},
  number={4},
  pages={617--639},
  year={2020},
  publisher={Springer}
}

@article{NW05,
  title={The algebra of {B}onferroni bounds: discrete tubes and extensions},
  author={Naiman, D. Q and Wynn, H. P},
  journal={Metrika},
  volume={62},
  number={2},
  pages={139--147},
  year={2005},
  publisher={Springer}
}

@article{PS21,
  title={Algebraic Analysis of Variants of Multi-State k-out-of-n Systems},
  author={Pascual-Ortigosa, P. and S{\'a}enz-de-Cabez{\'o}n, E.},
  journal={Mathematics},
  volume={9},
  number={17},
  pages={2042},
  year={2021},
  publisher={MDPI}
}

@article{PS22,
  title={Support posets of some monomial ideals},
  author={Pascual-Ortigosa, P. and {S{\'a}enz-de-Cabez{\'o}n}, E.},
  journal={Applicable Algebra in Engineering, Communication and Computing},
  volume={33},
  number={4},
  pages={457--475},
  year={2022},
  publisher={Springer}
}

@article{PSW21,
  title={Algebraic reliability of multi-state k-out-of-n systems},
  author={Pascual-Ortigosa, P. and S{\'a}enz-de-Cabez{\'o}n, E. and Wynn, H. P},
  journal={Probability in the Engineering and Informational Sciences},
  volume={35},
  number={4},
  pages={903--927},
  year={2021},
  publisher={Cambridge University Press}
}

@ARTICLE{P02,
	author = {Peeva, I.},
	title = {Resolutions and lattices},
	year = {2002},
	journal = {Homology, homotopy and applications},
	volume = {4},
	pages = {427--437}
}

@ARTICLE{P93,
	author = {Poincaré, H.},
	title = {Sur la généralisation d'un thèoreme d'{E}uler relatif aux polyèdres},
	year = {1893},
	journal = {C.R. Acad. Sci. Paris},
	volume = {117},
	pages = {144--145}
}

@book{R90,
  author    = {Ramamurthy, K. G.},
  title     = {Coherent Structures and Simple Games},
  publisher = {Kluwer Academic Publishers},
  address   = {Dordrecht},
  year      = {1990}
}

@article{R64,
  author  = {Rota, G.-C.},
  title   = {On the foundations of combinatorial theory {I}. {T}heory of {M}{\"o}bius functions},
  journal = {Zeitschrift f{\"u}r Wahrscheinlichkeitstheorie und Verwandte Gebiete},
  volume  = {2}, number = {4}, pages = {340--368}, year = {1964}
}

@article{RS13,
  title={Complexity and algorithms for Euler characteristic of simplicial complexes},
  author={Roune, Bjarke Hammersholt and S{\'a}enz-de-Cabez{\'o}n, Eduardo},
  journal={Journal of Symbolic Computation},
  volume={50},
  pages={170--196},
  year={2013},
  publisher={Elsevier}
}

@article{SP78,
  title={New topological formula and rapid algorithm for reliability analysis of complex networks},
  author={Satyanarayana, A. and Prabhakar, A.},
  journal={IEEE Transactions on Reliability},
  volume={27},
  number={2},
  pages={82--100},
  year={1978},
  publisher={IEEE}
}

@article{SC83,
  title={Network reliability and the factoring theorem},
  author={Satyanarayana, A. and Chang, M. K},
  journal={Networks},
  volume={13},
  number={1},
  pages={107--120},
  year={1983},
  publisher={Wiley Online Library}
}

@article{S09,
  title={Multigraded {B}etti numbers without computing minimal free resolutions},
  author={S{\'a}enz-de-Cabez{\'o}n, E.},
  journal={Applicable Algebra in Engineering, Communication and Computing},
  volume={20},
  number={5},
  pages={481--495},
  year={2009},
  publisher={Springer}
}

@article{SW09,
  title={Betti numbers and minimal free resolutions for multi-state system reliability bounds},
  author={S{\'a}enz-de-Cabez{\'o}n, E. and Wynn, H. P},
  journal={Journal of symbolic computation},
  volume={44},
  number={9},
  pages={1311--1325},
  year={2009},
  publisher={Elsevier}
}

@article{SW10,
  title={Mincut ideals of two-terminal networks},
  author={S{\'a}enz-de-Cabez{\'o}n, E. and Wynn, H. P},
  journal={Applicable Algebra in Engineering, Communication and Computing},
  volume={21},
  pages={443--457},
  year={2010},
  publisher={Springer}
}

@article{SW15,
  title={Hilbert functions in design for reliability},
  author={S{\'a}enz-de-Cabez{\'o}n, E. and Wynn, H. P},
  journal={IEEE transactions on reliability},
  volume={64},
  number={1},
  pages={83--93},
  year={2014},
  publisher={IEEE}
}

@article{SW85,
  author  = {Satyanarayana, A. and Wood, R. K.},
  title   = {A linear-time algorithm for computing {$K$}-terminal reliability in series-parallel networks},
  journal = {SIAM Journal on Computing},
  volume  = {14},
  number  = {4},
  pages   = {818--832},
  year    = {1985}
}

@article{S75,
  title={The upper bound conjecture and {C}ohen-{M}acaulay rings},
  author={Stanley, R. P},
  journal={Studies in Applied Mathematics},
  volume={54},
  number={2},
  pages={135--142},
  year={1975},
  publisher={Wiley Online Library}
}

@book{S96,
  title={Combinatorics and commutative algebra},
  author={Stanley, R. P},
  volume={41},
  year={1996},
  publisher={Birkh\"auser}
}

@article{S17,
  title={How the upper bound conjecture was proved},
  author={Stanley, R. P},
  journal={Selected Works of Richard P. Stanley},
  volume={25},
  pages={3},
  year={2017},
  publisher={American Mathematical Soc.}
}

@phdthesis{T66,
  title={Ideals generated by monomials in an R-sequence},
  author={Taylor, D. K.},
  year={1966},
  school={The University of Chicago}
}

@book{TB17,
  title={Reliability and availability engineering: modeling, analysis, and applications},
  author={Trivedi, K. S and Bobbio, A.},
  year={2017},
  publisher={Cambridge University Press}
}

@book{V15,
  title={Monomial algebras},
  author={Villarreal, R. H.},
  edition={2},
  year={2015},
  publisher={Chapman and Hall/CRC}
}

@article{W85,
  title={Multistate block diagrams and fault trees},
  author={Wood, A. P},
  journal={IEEE Transactions on Reliability},
  volume={34},
  number={3},
  pages={236--240},
  year={1985},
  publisher={IEEE}
}

@article{W86,
  author  = {Wood, R. K.},
  title   = {Factoring algorithms for computing {$K$}-terminal network reliability},
  journal = {IEEE Transactions on Reliability},
  volume  = {R-35},
  number  = {3},
  pages   = {269--278},
  year    = {1986}
}

@article{ZLW00,
  title={Reliability evaluation of combined k-out-of-n: F, consecutive-k-out-of-n: F and linear connected-(r, s)-out-of-(m, n): F system structures},
  author={Zuo, M. J. and Lin, D. and Wu, Y.},
  journal={IEEE Transactions on Reliability},
  volume={49},
  number={1},
  pages={99--104},
  year={2000},
  publisher={IEEE}
}

@article{ZT06,
  title={Performance evaluation of generalized multi-state k-out-of-n systems},
  author={Zuo, M. J. and Tian, Z.},
  journal={IEEE Transactions on Reliability},
  volume={55},
  number={2},
  pages={319--327},
  year={2006},
  publisher={IEEE}
}

\end{document}